\documentclass[11pt]{article}

\usepackage{amsmath, amsthm, amssymb, enumerate}
\usepackage{graphicx}

\date{}
\title{\vspace{-0.8cm}Bandwidth theorem for random graphs}
\author{
Hao Huang \thanks{Department of Mathematics, UCLA, Los
Angeles, CA, 90095. Email: huanghao@math.ucla.edu.}
\and
Choongbum Lee \thanks{Department of Mathematics, UCLA, Los
Angeles, CA, 90095. Email: choongbum.lee@gmail.com. Research
supported in part by Samsung Scholarship.}
\and
Benny Sudakov \thanks{Department of Mathematics, UCLA, Los Angeles, CA 90095.
Email: bsudakov@math.ucla.edu. Research supported in part by NSF
CAREER award DMS-0812005 and by a USA-Israel BSF grant. }
}

\theoremstyle{plain}
\newtheorem{thm}{Theorem}[section]
\newtheorem{prop}[thm]{Proposition}
\newtheorem{lemma}[thm]{Lemma}
\newtheorem{cor}[thm]{Corollary}
\newtheorem{claim}[thm]{Claim}

\theoremstyle{definition}
\newtheorem{dfn}{Definition}

\newenvironment{pf}{\noindent\textbf{Proof.}}{\qed\medskip}
\newenvironment{pfof}[1]{\noindent{\bf Proof of #1.}}{\qed\medskip}
\newenvironment{rem}{\noindent\textbf{Remark.}}{\medskip}

\newcommand{\floor}[1]{\left\lfloor #1 \right\rfloor}
\newcommand{\ceil}[1]{\left\lceil #1 \right\rceil}
\renewcommand{\deg}{\textrm{d}}

\begin{document}
\maketitle

\begin{abstract}
A graph $G$ is said to have \textit{bandwidth} at most $b$, if there
exists a labeling of the vertices by $1,2,\ldots, n$, so that $|i -
j| \leq b$ whenever $\{i,j\}$ is an edge of $G$. Recently,
B\"{o}ttcher, Schacht, and Taraz verified a conjecture of
Bollob\'{a}s and Koml\'{o}s which says that for every positive
$r,\Delta,\gamma$, there exists $\beta$ such that if $H$ is an
$n$-vertex $r$-chromatic graph with maximum degree at most $\Delta$
which has bandwidth at most $\beta n$, then any graph $G$ on $n$
vertices with minimum degree at least $(1 - 1/r + \gamma)n$ contains
a copy of $H$ for large enough $n$. In this paper, we extend this
theorem to dense random graphs. For bipartite $H$, this answers an open question of B\"{o}ttcher,
Kohayakawa, and Taraz. It appears that for non-bipartite $H$ the direct extension is not possible, and 
one needs in addition that some vertices of $H$ have independent neighborhoods.
We also obtain an asymptotically tight bound for the maximum number of vertex disjoint copies of a fixed
$r$-chromatic graph $H_0$ which one can find in a spanning subgraph
of $G(n,p)$ with minimum degree $(1-1/r + \gamma)np$.
\end{abstract}

\section{Introduction}
\label{section_introduction} One of the central themes in extremal
graph theory is the study of sufficient conditions which imply that
a graph $G$ contains a copy of a particular graph $H$. Two main
interesting cases of this problem are when $H$ has fixed order, and when
it has size comparable or the same as graph $G$. The celebrated
 Erd\H{o}s-Stone theorem \cite{MR0018807} settled the first case, showing that
sufficiently large graph $G$ of $n$ vertices and more than
$(1-\frac{1}{r-1}+o(1))\binom{n}{2}$ edges contains a copy of any
$r$-chromatic graph $H$ of fixed order.

In the second case, when the order of $H$ is close to the order of
$G$, the large number of edges is no longer sufficient to embed $H$
because there might be isolated vertices in $G$. Therefore we need a
lower bound on the minimum degree of $G$. The most well-known
example of such a result is Dirac's theorem (see, e.g.,
\cite{MR2159259}), which says that, if $G$ is a graph on $n$
vertices with minimum degree at least $\ceil{n/2}$ then $G$ contains
a Hamilton cycle. Another example is a problem of packing vertex
disjoint copies of a fixed graph $H_0$ in $G$. We say that $G$
contains a perfect $H_0$-packing if there are vertex disjoint copies
of $H_0$ that cover all the vertices of $G$. For convenience, we may
assume that the order of $G$ is divisible by the order of $H_0$. A
classical theorem of Hajnal and Szemer\'{e}di \cite{MR0297607}
states that if $G$ has minimum degree at least $(1- 1/r)n$ then $G$
contains a perfect packing of complete graphs $K_r$. More general
packing problems have been studied in \cite{MR1376050, KuOs2,
MR1829855}.

The $r$-th power of a graph $G$ is the graph $G^{(r)}$ obtained from
$G$ by connecting every pair of vertices which have distance at most
$r$ in $G$. In particular, note that the $(r-1)$-st power of the
$n$-cycle contains $\lfloor n/r \rfloor$ vertex disjoint copies of
$K_r$. P\'{o}sa and Seymour \cite{MR000005} proposed a common
generalization of Dirac's and Hajnal-Szemer\'{e}di's theorem. They
conjectured that the same minimum degree bound $(1-1/r)n$ will force
a graph $G$ to have the $(r-1)$-st power of a Hamiltonian cycle in
it. This conjecture has been open for quite a while until
Koml{\'o}s, S{\'a}rk\"{o}zy, and Szemer{\'e}di \cite{MR1682919}
proved it for large enough $n$. They used a combination of
Szemer\'{e}di regularity lemma \cite{MR540024} and the so-called
blow-up lemma \cite{MR1466579}. We will discuss this technique in
more detail later in the paper.

The above results might suggest that if $G$ has minimum degree at
least $(1 - 1/r + o(1))n$, then it contains a copy of any $n$-vertex
$r$-chromatic graph $H$ with bounded degree. However, the following
example (see, \cite{MR2448444}) shows that some restrictions are
necessary. Let $H$ be a random bipartite graph with bounded maximum
degree and parts of size $n/2$ and $G$ be a graph formed by two
cliques each of size $(1/2 + \gamma)n$ which share $2\gamma n$
vertices (for some small fixed $\gamma > 0$). Assume that $H$ is embedded into $G$ and look at the $(1/2
- \gamma)n$ vertices which come from one of the cliques and do not
belong to their intersection. The only neighbors of these vertices
in $G$ are the $2\gamma n$ vertices in the intersection. But with
high probability $H$ contains no collection of $(1/2 - \gamma)n$
vertices which have at most $2\gamma n$ neighbors. Therefore we
cannot embed $H$ into $G$.

Thus to find a general theorem, we need some additional restriction
on the graph $H$. A graph $H$ is said to have \textit{bandwidth} at
most $b$, if there exists a labeling of the vertices by $1,2,\ldots,
n$, so that $|i - j| \leq b$ whenever $i,j$ forms an edge. We denote
by $bw(H) = b$ if $b$ is the minimum integer such that $H$ has
bandwidth at most $b$. Bollob\'{a}s and Koml\'{o}s \cite{MR1684627}
conjectured that if $H$ is an $r$-chromatic graph which has bounded
degree and low enough bandwidth then one can embed it into a graph
$G$ with minimum degree at least $(1 - 1/r + o(1))n$. Note that the
constant $1-1/r$ is the best constant we can expect for such an
embedding result to hold. Indeed, assume that $n$ is divisible by
$r$ and let $G$ be the complete $r$-partite graph on $n$ vertices
whose partition classes are of size $n/r + 1, n/r - 1, n/r, \ldots,
n/r$. This graph has minimum degree $(1 - 1/r)n - 1$. Consider the
graph $H$ consisting of $n/r$ vertex disjoint copies of $K_r$. It is
clear that each copy of $K_r$ must contain at least one vertex from
each class of $G$ and thus there can only be at most $n/r - 1$ such
copies in $G$. Thus we cannot embed $H$ into $G$.

Bollob\'{a}s and Koml\'{o}s' conjecture has been recently proved by
B\"{o}ttcher, Schacht, and Taraz \cite{BoScTa2}, \cite{MR2448444}: for every positive
$r, \Delta, \gamma$, there exists $\beta$ such that if $H$ is an
$n$-vertex $r$-chromatic graph with maximum degree at most $\Delta$
and bandwidth at most $\beta n$, then any graph $G$ on $n$ vertices
with minimum degree at least $(1-1/r+\gamma)n$ contains a copy of $H$
for large enough $n$ (we will refer to this conjecture and theorem as the
bandwidth conjecture and the bandwidth theorem from now on). There
are a lot of graphs $H$ satisfying the condition above. For example, $r$-th powers of
cycles which have bandwidth $2r$, trees with constant maximum degree
which have bandwidth at most $O(n / \log n)$ \cite{MR1205400}, and
$n^{1/2}$ by $n^{1/2}$ square grids which have bandwidth $O(n^{1/2})$
are a few of those. For more examples, see
B\"{o}ttcher, Pruessmann, Taraz, and W\"{u}rfl's \cite{MR000003}
classification of bounded degree graphs with sublinear bandwidth.
Moreover, the theorem proved in \cite{MR2448444} is a strengthening
of the bandwidth conjecture and also implies Dirac's theorem and
P\'{o}sa-Seymour's conjecture asymptotically.

Most of the above mentioned results can also be viewed in the framework
of \textit{resilience} which we discuss next.
 A graph property is called \textit{monotone
increasing (decreasing)} if it is preserved under edge addition
(deletion). Following \cite{MR2462249}, we define:

\begin{dfn}
Let $\mathcal{P}$ be a monotone increasing (decreasing) graph
property.\\
\begin{tabular}{cp{15.5cm}}
(i) & The global resilience of $G$ with respect to $\mathcal{P}$ is the minimum number $r$ such that by deleting (adding) $r$ edges from $G$ one can obtain a graph not having $\mathcal{P}$.\\
(ii) & The local resilience of a graph $G$ with respect to $\mathcal{P}$ is the minimum number $r$ such that by deleting (adding) at most $r$ edges at each vertex of $G$ one can obtain a graph not having $\mathcal{P}$.\\
\end{tabular}
\end{dfn}

Intuitively, the question of determining resilience of a graph $G$
with respect to a graph property $\mathcal{P}$ is like asking, ``How
strongly does $G$ possess $\mathcal{P}$?''. Using this terminology,
one can for example restate Dirac's theorem as saying that $K_n$ has
local resilience $\lfloor n/2 \rfloor$ with respect to having a
Hamilton cycle. In \cite{MR2462249}, Sudakov and Vu have initiated
the systematic study of global and local resilience of random and
pseudorandom graphs. The random graph model they considered is the
binomial random graph $G(n,p)$, which denotes the probability space
whose points are graphs with vertex set $[n] = \{1,\ldots,n\}$ where
each pair of vertices forms an edge randomly and independently with
probability $p$. Given a graph property $\mathcal{P}$, we say that
$G(n,p)$ possesses $\mathcal{P}$ \textit{asymptotically almost
surely}, or a.a.s. for brevity, if the probability that $G(n,p)$
possesses $\mathcal{P}$ tends to 1 as $n$ tends to infinity. In the
above mentioned paper, Sudakov
and Vu studied the resilience of random graphs with respect to various
properties such as Hamiltonicity, containing a perfect matching,
increasing its chromatic number,
and having a nontrivial automorphism (this result appeared in their
earlier paper with Kim \cite{MR1945368}). For example, they proved
that if $p > \log^4 n / n$ then a.a.s. any subgraph of $G(n,p)$ with
minimum degree $(1/2+o(1))np$ is Hamiltonian. Note that this result
can be viewed as a generalization of Dirac's theorem mentioned
above, since the complete graph is also a random graph $G(n,p)$ with
$p=1$. This connection is very natural and most of the resilience
results for random and pseudorandom graphs can be viewed as a
generalization of classical results from graph theory. For
additional resilience type results, see, e.g. \cite{MR000007,
MR000000, MR000006,MR2383452, MR2430433, MR000001}.

Using the above terminology, the bandwidth theorem says that the
complete graph $K_n$ has local resilience $(1/r+o(1))n$ with respect
to containing spanning $r$-chromatic graphs $H$ of low bandwidth and bounded
degree. B\"{o}ttcher, Kohayakawa, and Taraz \cite{MR000006}
partially extended this result to random graphs by proving that for
fixed $\eta, \gamma > 0, \Delta > 1$ there exist positive constants
$\beta$ and $c$ such that if $p \geq c(\log n / n)^{1 / \Delta}$
then a.a.s every subgraph of $G(n,p)$ with minimum degree at least
$(1/2 + \gamma)np$ contains a copy of any bipartite graph $H$ with
$(1-\eta)n$ vertices, maximum degree $\Delta$ and bandwidth at most
$\beta n$. They then posed a natural
and interesting question \cite{MR000004}, whether one can
fully extend the bandwidth theorem to random graphs. More
specifically, they suggested that it should be possible to extend the
bandwidth theorem for spanning bipartite $H$ in the regime of constant edge probability $p$.
For this range of probabilities, there are well developed tools that
we can use, and thus there are more hopes to understand the correct
behavior of this problem. The reason we only focus on bipartite
graphs is the following. Consider the problem of finding a triangle
factor. A fixed vertex $v$ in $G(n,p)$ a.a.s.~has degree
$(1+o(1))np$ and has $(1+o(1))np^2$ common neighbors with any other
vertex. If we delete all the edges in the neighborhood of $v$, we
destroy all the triangles containing $v$. On the other hand, the
degree of any vertex in $G(n,p)$ will decrease by at most
$O(np^2)\ll np$, and thus, it will still be greater than
$(2/3+\gamma)np$. This gives a subgraph of $G(n,p)$ with minimum degree at least 
$(2/3+\gamma)np$ and no triangle factor. 
Since disjoint union of triangles has constant bandwidth, this simple observation shows that one can not directly 
extend the bandwidth theorem in full generality.

In this paper we study the above mentioned question posed by B\"{o}ttcher,
Kohayakawa, and Taraz. We have the following two main contributions. First, we prove
that for constant edge probability, it is possible to obtain a complete
extension of the bandwidth theorem for spanning
bipartite graphs $H$ with bounded degree and sublinear bandwidth.
We also suggest a natural minor restriction on  
non-bipartite graphs $H$, which makes possible an extension of bandwidth theorem to random graphs. 
More precisely, we show that having some vertices with independent neighborhoods in $H$ is enough.
Here our main theorem.

\begin{thm}\label{thm_mainthmintro1}
For fixed integers $r, \Delta$, and reals $0 < p \leq 1$ and $\gamma
> 0$, there exists a constant $\beta > 0$ such that a.a.s., any
spanning subgraph $G'$ of $G(n,p)$ with minimum degree $\delta(G')
\ge (1 - 1/r + \gamma)np$ contains every $n$-vertex graph $H$ which
satisfies the following properties. $(i)$ $H$ is $r$-chromatic,
$(ii)$ has maximum degree at most $\Delta$, $(iii)$ has bandwidth at
most $\beta n$ with respect to a labeling of vertices by
$1,2,\ldots, n$, and $(iv)$ for every interval $[a, a+\beta^2 n]
\subset [1,n]$, there exists a vertex $v \in H$ such that $N_H(v)$
is an independent set.
\end{thm}

In particular, the theorem holds for any bipartite $H$ which has
bounded degree and sublinear bandwidth. Thus it positively answers the 
above mentioned question of B\"{o}ttcher, Kohayakawa, and
Taraz for dense random graphs. Note that for non-bipartite
graphs, we only require constant number of vertices 
with independent neighborhoods. 

Another main contribution of this paper is an extension of the classical extremal results
on $H_0$-packings in graphs with large minimum degree to the setting of 
random graphs. The above theorem implies that if $H_0$ is a fixed $r$-chromatic graph
having a vertex not contained in a triangle and $\gamma>0$ is any
fixed constant, then a.a.s every $G'\subset G(n,p)$ with minimum
degree at least $(1-1/r+\gamma)np$ contains a perfect $H_0$-packing.
This suggests the following natural question. Let $H_0$ be a fixed
$r$-chromatic graph whose every vertex belongs to some triangle.
What is the maximum number of vertex disjoint copies of $H_0$ that
one can find in a spanning subgraph $G'$ of $G(n,p)$ with
$\delta(G') \geq (1-1/r + \gamma)np$? We proved the following
result, which gives a rather accurate answer to this question.

\begin{thm} \label{thm_packingtheorem}
Let $H_0$ be an $r$-chromatic graph whose every vertex is contained
in a triangle. Then there exist constants $c=c(r)$ and $C=C(r)$
such that for any fixed $0 < p \leq 1$ and $0 < \gamma \leq 1/(2r)$, the
random graph $G(n,p)$ a.a.s.~has the following properties.

$(i)$ There exists a spanning subgraph $G'$ with minimum degree
$\delta(G') \geq (1 - 1/r + \gamma)np$ such that at least
$\left\lfloor cp^{-2} \right\rfloor$ vertices of $G'$ are not
contained in a copy of $H_0$.

$(ii)$ For every spanning subgraph $G'\subset G$ which has minimum
degree $\delta(G')\geq (1-1/r+\gamma)np$, at least $n - Cp^{-2}$
vertices of $G'$ can be covered by vertex disjoint copies of $H_0$.
\end{thm}

The rest of this paper is organized as follows. In Section
\ref{section_preliminaries} we collect some known results which we
need later to prove our main theorem. In Section
\ref{section_outlineproof} we state several important lemmas, and
outline the proof of the main theorem using these lemmas. In Section
\ref{section_thekeylemma} we prove the lemmas given in
Section \ref{section_outlineproof}. In Section \ref{section_maintheorem}
we provide a detailed proof of Theorem \ref{thm_mainthmintro1} by
using the tools developed in previous sections and other known
results. As an application of the main theorem, in Section
\ref{section_applicationpackingproblem}, we study the packing
problem in random graphs (Theorem \ref{thm_packingtheorem}). The
last section contains some concluding remarks and open problems. In
the appendix, we extend the main theorem to pseudorandom graphs
(only the sketch of the proof will be given).

To simplify the presentation, we often omit floor and ceiling signs
whenever these are not crucial and make no attempts to optimize
absolute constants involved. We also assume that the order $n$ of
all graphs tends to infinity and therefore is sufficiently large
whenever necessary. Throughout the paper, whenever we refer, for
example, to a function with
subscript as $f_{3.1}$, we mean the function $f$ defined in Lemma/Theorem 3.1.\\

\noindent\textbf{Notation.} $G=(V,E)$ denotes a graph with vertex
set $V$ and edge set $E$. $\Delta(G), \delta(G), \chi(G)$ denote the
maximum degree, the minimum degree, and the chromatic number of $G$
respectively. In the following, we will use $v$ for a vertex and $X$
for an arbitrary set. Let $N(X)$ be the collection of all vertices
which are adjacent to at least one vertex in $X$. If $X=\{v\}$ is a
singleton set we denote its neighborhood by $N(v)$. Let $N^{(0)}(v)
:= \{ v \}$ and $N^{(k)}(v)$ be the vertices at distance exactly $k$
from $v$. Note that $N^{(1)}(v) = N(v)$. Similarly define
$N^{(k)}(X)$ to be the vertices at distance exactly $k$ from the set
$X$, where the distance of a vertex $v$ from a set $X$ is defined as
the minimum number $t$ such that $N^{(t)}(v) \cap X \neq \emptyset$.
The degree of a vertex is defined as $\deg(v) := |N(v)|$. The
neighborhood of a vertex in a set is defined as $N(v,X) := N(v) \cap
X$ and the degree of a vertex in a set is defined as $\deg(v,X) :=
|N(v,X)|$. We denote by $E(X)$ the set of edges in the induced
subgraph $G[X]$ and by $e(X):= |E(X)|$ its size. Similarly, for two
sets $X$ and $Y$, we denote by $E(X,Y)$ the set of ordered pairs
$(x,y) \in E$ such that $x \in X$ and $y \in Y$, also $e(X,Y) :=
|E(X,Y)|$. Note that $e(X,X) = 2e(X)$. By $d(X,Y) := e(X,Y)/|X||Y|$
we denote the density of the pair. If we have several graphs, then
the graph we are currently working with will be stated as a
subscript. For example $N^{(k)}_G(v)$ is the $k$-th neighborhood of
$v$ in graph $G$.

We also utilize the following standard asymptotic notation. For two
functions $f(n)$ and $g(n)$, write $f(n)=\Omega(g(n))$ if there
exists a constant $C$ such that $\liminf_{n \rightarrow \infty}
f(n)/g(n) \geq C$. If there is a subscript such as in
$\Omega_\varepsilon$ this means that the constant $C$ may depend on
$\varepsilon$. We write $f(n)=o(g(n))$ or $f(n) \ll g(n)$  if
$\limsup_{n\rightarrow \infty}f(n)/g(n) = 0$. Also, $f(n)=O(g(n))$
if there exists a positive constant $C>0$ such that
$\limsup_{n\rightarrow \infty} f(n)/g(n) \leq C$. Throughout the
paper log denotes the natural logarithm.

\section{Preliminaries}
\label{section_preliminaries}
In this section, we collect several known results to be used later
in the proof of the main theorem.

The following well-known concentration result (see, for example
\cite{AlSp}, Appendix A) will be used several times throughout the
proof. We denote by $Bi(n,p)$ a binomial random variable with
parameters $n$ and $p$.
\begin{thm}[Chernoff Inequality] \label{thm_chernoffinequality}
If $X \sim Bi(n,p)$ and $\lambda \leq np$, then
\[ P \big( |X - np| \geq \lambda \big) \leq e^{-\Omega(\lambda^2/(np))}. \]
\end{thm}

Our approach in proving the main theorem is to use the regularity
lemma and the blow-up lemma. These powerful tools developed by
Szemer{\'e}di \cite{MR540024}, and Koml{\'o}s, S{\'a}rk\"{o}zy,
Szemer{\'e}di \cite{MR1466579}, respectively, have been successively
applied to solve several embedding results (e.g.,
\cite{KuOs}). Here we state
these facts without proof. Readers may consult \cite{MR1395865},
\cite{MR1684627} for more detailed discussion on these topics.

Let $G=(V,E)$ be a graph and $\varepsilon > 0$ be fixed. A disjoint
pair of sets $X,Y \subset V$ is called an {\em $\varepsilon$-regular
pair in $G$} if all $A \subset X, B \subset Y$ such that $|A| \ge
\varepsilon |X|, |B| \ge \varepsilon |Y|$ satisfy $|d(X,Y) - d(A,B)|
\leq \varepsilon$. An $\varepsilon$-regular pair $(X,Y)$ is called
$(d, \varepsilon)$-regular, if it has density at least $d$. A vertex
partition $V_0, \ldots, V_k$ is called an {\em $\varepsilon$-regular
partition of $G$} if (i) $|V_0| \leq \varepsilon n$, (ii) $V_i$ have
equal size for $i \geq 1$, and (iii) $(V_i, V_j)$ is
$\varepsilon$-regular in $G$ for all but at most $\varepsilon k^2$
pairs $1 \leq i < j \leq n$.  The regularity lemma states that every
large enough graph admits a regular partition. Here we state it in a
stronger form which can be found in \cite{MR1395865}.

\begin{lemma}[Regularity Lemma] \label{thm_regularitylemma}
For every integer $t$ and real $\varepsilon > 0$, there exists $n_0 = n_0(t,
\varepsilon)$ and $T = T(t,\varepsilon)$ such that for every graph $G$ on $n \geq n_0$
vertices and $d \in [0,1]$, there exists a subgraph $G' \subset G$
with an $\varepsilon$-regular partition $V_0, \ldots, V_k$ of $G'$
satisfying the following properties.\\
\begin{tabular}{cl}
$(i)$ & $t \le k \le T$,\\
$(ii)$ & $\deg_{G'}(v) > \deg_G(v) - (d + \varepsilon)n$ for all $v \in V$.\\
$(iii)$ & $e(G'[V_i]) = 0$ for all $i \geq 1$,\\
$(iv)$ & every pair $(V_i, V_j) \,(1 \leq i < j \leq k)$ either is $\varepsilon$-regular in $G'$ with density at least $d$ or has \\
& no edges between them.
\end{tabular}
\end{lemma}

Let $V_0, \ldots, V_k$ be an $\varepsilon$-regular partition of $G$.
Then we define the {\em reduced graph} $R$ with parameters $(d,
\varepsilon)$ as the graph on the vertex set $[k]$ with edges
$\{i,j\} \in E(R)$ if and only if $(V_i, V_j)$ is
$(d,\varepsilon)$-regular. In this case, we also say that $V_0,
\ldots, V_k$ is {\em $(d, \varepsilon)$-regular on $R$ in $G$}.
Furthermore, if $G' \subset G$ and $V_0, \ldots, V_k$ satisfy
$(i),(ii),(iii),(iv)$ of Lemma \ref{thm_regularitylemma}, we say
that $V_0, \ldots, V_k$ is a {\em pure $(d,\varepsilon)$-regular
partition of $G'$}. The following lemma establishes the fact that
the reduced graph inherits the minimum degree condition.

\begin{lemma} \label{lemma_reducedgraphmindegree}
Let $0< p \le 1$ and $\alpha, \gamma>0$ be fixed. There exists
$\epsilon_0=\epsilon_0(p,\alpha,\gamma)$ such that for all
$\varepsilon \le \varepsilon_0$ and $d > 0$, the following
a.a.s.~holds. Given a graph $G=G(n,p)$, let $V_0, V_1, \ldots, V_k$
be a pure $(d,\varepsilon)$-regular partition of a subgraph $G'
\subset G$, and $R$ be its reduced graph. If $G'$ has minimum degree
at least $(\alpha + \gamma)np$, then $R$ has minimum degree at least
$(\alpha + 3\gamma/4)k$.
\end{lemma}
\begin{pf}
Let $m := |V_i|$. Since $|V_0| \le \varepsilon n$, we have the bound
$m \geq (1-\varepsilon)n/k$. Thus by Chernoff inequality,
a.a.s.~$e_{G'}(V_i, V_j) \leq e_G(V_i, V_j) \leq (1+
\varepsilon)m^2p$ for all $i,j \ge 1$. From the definition of a pure
$(d,\varepsilon)$-regular partition, we know that if $\{i,j\} \notin
E(R)$, then the pair $(V_i,V_j)$ has no edges between them. Thus for
a vertex $i \in V(R)$,
\[
 e_{G'}(V_i, V \setminus V_0) = \sum_{j=1}^{k} e_{G'}(V_i, V_j) \leq (k - \deg_R(i))\cdot 0 + \deg_R(i)(1+\varepsilon)m^2p = \deg_R(i)(1+\varepsilon)m^2p.
\]
On the other hand, by the minimum degree condition of $G'$ and the fact
$e_{G'}(V_i) = 0$,
\[ e_{G'}(V_i, V
\setminus V_0) \ge \left(\sum_{v \in V_i}\deg_{G'}(v) \right) -
e_{G'}(V_i, V_0) \ge (\alpha + \gamma)np|V_i| - \varepsilon n |V_i| .
\]
Combine the bounds, divide each side by $m^2p$ and use the bound $n >
mk$ to get, $(\alpha + \gamma - \varepsilon/p )k \le (1+\varepsilon)\deg_R(i)$.
By selecting $\varepsilon$ small
enough, we have $\deg_R(i) \ge (\alpha + 3\gamma/4)k$.
\end{pf}

With respect to embedding small subgraphs, regular pairs behave like
random graphs. Thus, merely knowing the structure of the reduced
graph already tells us plenty of information about the original
graph and the subgraphs that it contain. The following lemma is a
formal description of this intuition. A {\em graph homomorphism}
between two graphs $G_1 = (V_1, E_1), G_2 = (V_2, E_2)$ is a map $f:
V_1 \rightarrow V_2$ such that $(f(v),f(w)) \in E_2$ if $(v,w) \in
E_1$. We say that $G_1$ is homomorphic to $G_2$ if there is a
homomorphism from $G_1$ to $G_2$.

\begin{thm} \label{thm_embeddinglemma}
For any fixed graph $H$ and $d>0$, there exists an $\varepsilon_0 >
0$ such that for all $\varepsilon \le \varepsilon_0$, there is an
$n_0$ with the following property. Let $G$ be a graph on $n \ge n_0$
vertices, $V_0, \ldots, V_k$ be an $\varepsilon$-regular partition
of $G$, and $R$ be its reduced graph with parameters $(d,
\varepsilon)$. If $H$ is homomorphic to $R$, then $G$ contains a copy
of $H$.
\end{thm}

It is well known that the regularity lemma together with this
embedding lemma implies the following generalization of
Erd\H{o}s-Stone theorem to random graphs $G(n,p)$ when $p \in (0,1]$
is fixed (see, e.g., \cite{MR2187740} for discussion of the case $p
\ll 1$.) Recently, by using a different approach, Conlon and Gowers
\cite{CoGo}, and Schacht \cite{Schacht} independently extended this
result to the range $p \ll 1$, but we do not need this stronger form
for our purpose.

\begin{thm}  \label{thm_randomturan}
For any fixed $\gamma>0, 0 < p \leq 1$ and a graph $H$, $G=G(n,p)$
satisfies the following with probability $1 - e^{-\Omega(n^2p)}$.
Any subgraph $G' \subset G$ with
$e(G') \geq (1 - 1/(\chi(H)-1) + \gamma)n^2p/2$ contains a copy of
$H$.
\end{thm}

In fact, we need the following seemingly stronger result which
directly follows from Theorem \ref{thm_randomturan} by taking the
union bound.

\begin{cor}  \label{cor_randomturanstrong}
For any fixed $\alpha, \gamma>0, 0 < p \leq 1$ and a graph $H$,
$G=G(n,p)$ satisfies the following with probability $1 -
e^{-\Omega(n^2p)}$. For any subset $W \subset V$ of size $|W| \ge
\alpha n$, every subgraph $G' \subset G[W]$ with $e(G') \geq (1 -
1/(\chi(H)-1) + \gamma)|W|^2p/2$ contains a copy of $H$.
\end{cor}

The theorems above illustrate the strength of regularity in finding
fixed size subgraphs. On the other hand, the blow-up lemma, which we
will introduce next, exemplifies the strength of regularity in
embedding graphs which are as large as $G$. Before we state the
theorem we must define the concept of {\em super-regularity}. Let
$G=(V,E)$ be a graph and $d,\varepsilon > 0$. Then a pair of
disjoint sets $X,Y \subset V$ is called {\em $(d,
\varepsilon)$-super-regular in $G$} if it is (i)
$(d,\varepsilon)$-regular in $G$, and (ii) $\forall x \in X,
\deg_Y(x) \geq d|Y|$ and $\forall y \in Y, \deg_X(y) \geq d|X|$. As
in the regularity case, given a partition $V_0, \ldots, V_k$ of $G$
we define the $(d,\varepsilon)$-super-regular reduced graph $R$ to
be the graph on the vertex set $[k]$ with edges $\{i,j\} \in E(R)$
if and only if $(V_i, V_j)$ forms a $(d,\varepsilon)$-super-regular
pair in $G$. We may also say that $V_0, \ldots, V_k$ is
$(d,\varepsilon)$-super-regular on $R$ in $G$. The following version
of the blow-up lemma was used in \cite{BoScTa2} and
\cite{MR2448444}.

\begin{thm}[Blow-up lemma] \label{thm_blowuplemma}
For any positive $d, \Delta, c$ and $r$, there exist $\varepsilon =
\varepsilon(d,\Delta,c,r)$ and $\alpha=\alpha(d,\Delta,c,r)$ such
that the following is true. Let $n_1, n_2, \cdots ,n_r$ be arbitrary
integers and consider the following two graphs over the vertex set
$V=V_1 \cup \cdots \cup V_r$ with $|V_i| = n_i$ for all $1 \leq i
\leq r$. $(i)$ In $G_0$, each pair $(V_i, V_j)$ forms a complete
bipartite graph, and $(ii)$ in $G_1$, each pair $(V_i,V_j)$ forms a
$(d, \varepsilon)$ super-regular pair. Then any graph $H=(W_1 \cup
\ldots \cup W_r, E_H)$ with $\Delta(H) \leq \Delta$ and $|W_i|=n_i$
($\forall i \in [r]$) which can be embedded into $G_0$ so that all
the vertices of $W_i$ get mapped into $V_i$ ($\forall i \in [r]$)
can be embedded into $G_1$ in the same way.

Moreover, assume that we are given subsets $W_i' \subset W_i$ such
that $|W_i'| \leq \alpha \cdot \min_{j \in [r]} |W_j|$, and for each
$w \in W_i'$, a set $C_{w} \subset V_i$ such that $|C_{w}| \geq c
|V_i|$. Then there exists an embedding of $H$ into $G$ such that
every vertex $w \in W_1' \cup \ldots \cup W_k'$ is mapped into a
vertex in $C_w$.
\end{thm}

\section{Outline of the proof}
\label{section_outlineproof}

The setting of Theorem \ref{thm_mainthmintro1} can be briefly stated
as following. We have a host graph $G' \subset G(n,p)$ with large
minimum degree, a graph $H$ with certain restrictions, and we want
to embed $H$ into $G'$. Hence, with this setting in mind, in the
future discussion, $G'$ will always stand for the host graph, and
$H$ will stand for the graph that we want to embed.

To prove Theorem \ref{thm_mainthmintro1}, we adapt several lemmas
from the proof of the bandwidth theorem given in \cite{MR2448444}.
In this section, we will provide the statement of the lemmas, and
outline the proof of the main theorem by using these lemmas. The
statement of these lemmas might seem quite technical, so to
understand the intuition which lies behind the lemmas, it will be
useful to keep in mind that the final part of the proof will be an
application of the blow-up lemma given in Theorem
\ref{thm_blowuplemma}.

\begin{figure}[t]
\centering
\includegraphics[height=2.0cm]{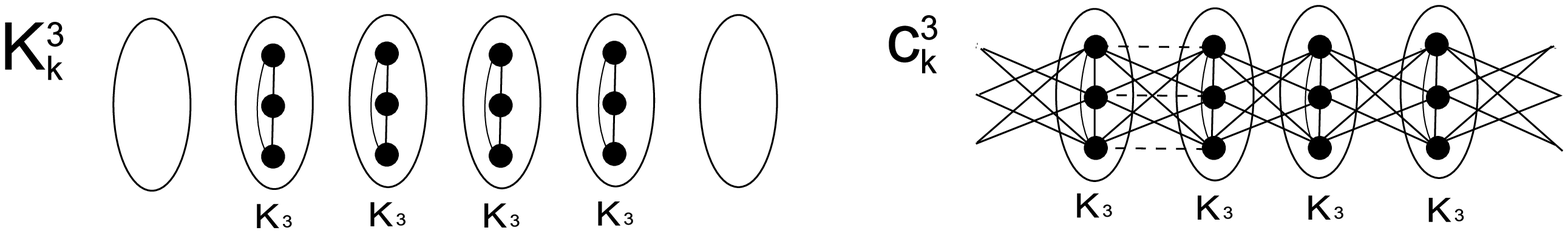}
\caption{$K^3_k$ and $C^3_k$}
 \label{figure_k3kc3k}
\end{figure}

First lemma, which is a variant of `Lemma for G'(Lemma 6 in
\cite{MR2448444}), prepares the graph $G'$ so that we have many
regular and super-regular pairs. Before stating the lemma, we
introduce some graphs. The graphs $C_k^r$ and $K_k^r$ are defined as
following (see figure \ref{figure_k3kc3k}). $C_k^r$ is a graph over
the vertex set $[k] \times [r]$ such that $(i_1, j_1), (i_2, j_2)$
is connected by an edge if (i) $i_1 = i_2$ and $j_1 \neq j_2$, or
(ii) $|i_2 - i_1| = 1$ and $j_1 \neq j_2$. $K_k^r$ is a graph over
the same vertex set $[k] \times [r]$ consisting of $k$ disjoint
copies of $K_r$ each of which lies on the vertices $\{ i \} \times
[r]$. Note that $K_k^r \subset C_k^r$ by construction.

An integer partition $(n_{i,j})_{1 \leq i \leq k, 1 \leq j \leq r}$
of $n$ is called $r$-equitable if $|n_{i,j} - n_{i, j'}| \le 1$
for all $1 \le i \le k$ and $1 \le j,j' \le r$.

\begin{lemma}[Lemma for G] \label{lemma_lemmaforG}
For every integer $r \ge 2$, $0 < p \leq 1$ and $\gamma > 0$ there exists
$d = d(r,p,\gamma) > 0$ and $\varepsilon_0=\varepsilon_0(r,p,\gamma)> 0$ such that for every
positive $\varepsilon \leq \varepsilon_0$ there exists
$b_0 = b_0(r,p,\gamma,\varepsilon)$,  $\xi_0 = \xi_0(r,p,\gamma,\varepsilon)
> 0$, and $K_0=K_0(r,p,\gamma,\varepsilon)$ such that, $G=G(n,p)$ a.a.s. satisfies the
following. For every subgraph $G' \subset G$ with $\delta(G') \geq
(1-1/r + \gamma)np$ there exist a subgraph $G'' \subset G'$ with
$\delta(G'') \geq (1-1/r + 4\gamma/5)np$, a set $B$ of size at most
$b_0$, an $r$-equitable integer partition $(m_{i,j})_{1 \leq i \leq k, 1 \leq j \leq r}$ of $n - |B|$,
sets $( V_{i,j}^* )_{1 \leq i \leq k, 1 \leq j \leq r}$,
and a graph $R$ on vertex set $[k] \times [r]$ with $k \leq K_0$ such that \\
\begin{tabular}{cl}
$(i)$ & $K_k^r \subset C_k^r \subset R$ and $\delta(R) \geq (1-1/r + \gamma/2)kr$, \\
$(ii)$ & $\forall 1 \le i \le k, 1 \le j \le r$, $m_{i,j} \geq (1-\varepsilon)n/(kr)$,\\
$(iii)$ & $\forall 1 \leq i \leq k, 1 \leq j \leq r$, $m_{i,j} \ge |V_{i,j}^*| \geq (1-\varepsilon)m_{i,j}$,\\
$(iv)$ & $(V_{i,j}^*)_{1 \leq i \leq k, 1 \leq j \leq r}$ is $(d,\varepsilon)$-regular on $R$ in $G''$, such that
\end{tabular}

\noindent for every choice of $(n_{i,j})_{1 \leq i \leq k, 1 \leq j
\leq r}$
 with $m_{i,j} - \xi_0 n \leq n_{i,j} \leq m_{i,j} + \xi_0 n$ and $\sum_{i,j} n_{i,j} \leq n - |B|$,
there exists a partition $(V_{i,j})_{1 \leq i \leq k, 1 \leq j \leq
r}$ of $V \backslash B$
with \\
\begin{tabular}{cl}
$(a)$ & $|V_{i,j}| \geq n_{i,j}$, $V_{i,j}^* \subset V_{i,j}$, $\forall 1 \leq i \leq k, 1 \leq j \leq r$, \\
$(b)$ & $(V_{i,j})_{1 \leq i \leq k, 1 \leq j \leq r}$ is $(d,\varepsilon)$-regular on $R$ in $G''$ and \\
$(c)$ & $(V_{i,j})_{1 \leq i \leq k, 1 \leq j \leq r}$ is $(d,\varepsilon)$-super-regular on $K_k^r$ in $G''$.
\end{tabular}
\end{lemma}

\begin{figure}[t]
\centering
\includegraphics[width=14cm]{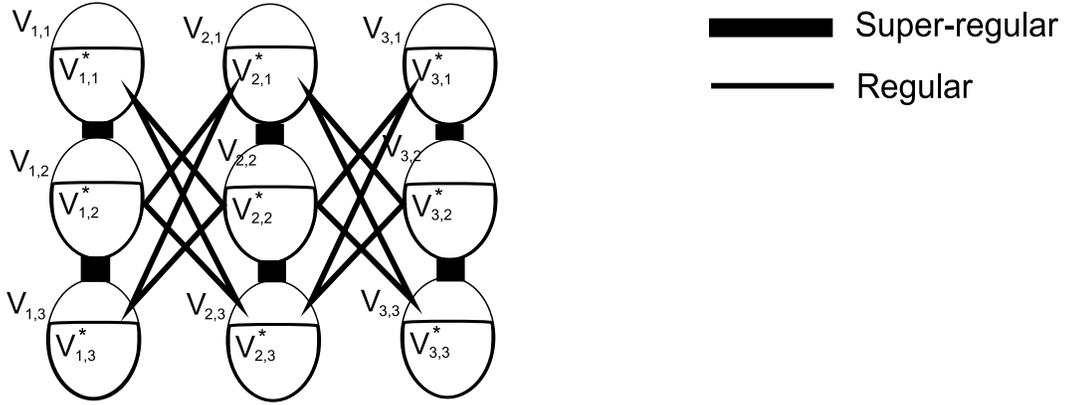}
\caption{Lemma for G. $V_{i,j}^*$ are fixed and $V_{i,j}$ are flexible in size.}
\label{figure_coreset}
\end{figure}

Heuristically, given a graph $G'$, this lemma returns some set $B$
and a `temporary' vertex partition of $V \backslash B$ with parts of
size $m_{i,j}$ for some integer partition $(m_{i,j})_{1 \leq i \leq
k, 1 \leq j \leq r}$ of $n - |B|$. The vertex partition is flexible
in the sense that given any other integer partition $(n_{i,j})_{1
\leq i \leq k, 1 \leq j \leq r}$ which is close to
$(m_{i,j})_{i,j}$, we can change the partition slightly so that the
new partition $(V_{i,j})_{i,j}$ has size $|V_{i,j}| = n_{i,j}$ for
all $i,j$. Moreover, each partition $V_{i,j}$ has an underlying
`core' set $V_{i,j}^*$ which always remains where they were
regardless of the given $(n_{i,j})_{i,j}$ (see figure
\ref{figure_coreset}). The main difference between this lemma and
`Lemma for G' in \cite{MR2448444} is the set $B$ whose existence is
unavoidable due to the inherent randomness of $G'$, and the `core'
sets $V_{i,j}^*$ which are there to help controlling the set $B$.
Note that $|B|$ is bounded by some constant $b_0$ which does not
depend on $n$.

Assume for the sake of argument, that the graph $H$ which we want to
embed into $G'$ consists of vertex disjoint copies of $C_4$, and
$r=2$, and $n$ is divisible by 4. Provide the graph $G'$ to Lemma
\ref{lemma_lemmaforG}, and get as output an integer partition $(m_{i,j})_{1
\leq i \leq k, 1 \leq j \leq 2}$ and a set $B$. If $B$ were empty,
then the rest of the argument can go as following. Find an integer
partition $(n_{i,j})_{1 \leq i \leq k, 1 \leq j \leq 2}$ which is
close to $(m_{i,j})_{1 \leq i \leq k, 1 \leq j \leq 2}$, and
satisfies $n_{i,1} = n_{i,2}$ with both $n_{i,1}, n_{i,2}$ being an even
integer for all $1 \leq i \leq k$. By Lemma \ref{lemma_lemmaforG},
we can obtain a partition $(V_{i,j})_{1 \leq i \leq k, 1 \leq j \leq
2}$ of $V \backslash B = V$ such that $|V_{i,j}| = n_{i,j}$ for all
$i,j$. Then apply the blow-up lemma on each copy of $K_2$ in $K_k^2$
separately, to find vertex disjoint copies
of $C_4$ in the graph.

To cover the case when $B$ is not empty, we need to slightly modify
this argument. As a first step, find copies of $C_4$ which only use
vertices from $B$ and $(V_{i,j}^*)_{1 \leq i \leq k, 1 \leq j \leq
2}$. Assume that there are no remaining vertices in $B$ after
finding some copies of $C_4$ (this part is not trivial but assume
that we can do this), and by doing so we have used $\delta_{i,j}$
vertices from each set $V_{i,j}^*$. Then $n - \sum_{i,j}
\delta_{i,j} - |B|$ is divisible by 4 and hence we can find an
integer partition $(n_{i,j})_{1 \leq i \leq k, 1 \leq j \leq 2}$ of
it which satisfies $n_{i,1} = n_{i,2}$ with both of
$n_{i,1},n_{i,2}$ being an even integer for all $1 \leq i \leq k$.
If this integer partition were also close to $(m_{i,j})_{1 \leq i
\leq k, 1 \leq j \leq 2}$, then by Lemma \ref{lemma_lemmaforG}, we
can obtain a partition $(V_{i,j})_{1 \leq i \leq k, 1 \leq j \leq
2}$ of $V \backslash B$ such that $|V_{i,j}| = n_{i,j} +
\delta_{i,j}$ for all $i,j$. Recall that the copies of $C_4$ which
we have already found use $\delta_{i,j}$ vertices from each set
$V_{i,j}^*$, and thus also from $V_{i,j}$. Therefore the remaining
number of vertices in $V_{i,j}$ after disregarding these copies of
$C_4$ is exactly $n_{i,j}$. Also note that deleting constant number
of vertices from each part does not destroy super-regularity. Now
apply the blow-up lemma to the remaining vertices and find vertex
disjoint copies of $C_4$ which cover all the vertices of $V$.

The strategy of embedding a general graph $H$ is not too different
from this. However, a general graph $H$ can have more complicated structure than
vertex disjoint copies of $C_4$, and requires some preprocessing
before being embedded into $G'$. In the next lemma, we use the bound
on the bandwidth to map $H$ `nicely' onto the $[k] \times [r]$ grid.
This lemma is a variant of `Lemma for H' (Lemma 8 in \cite{MR2448444}),
and can be derived from it without much difficulty.

\begin{lemma}[Lemma for H] \label{lemma_lemmaforH}
Let $r,k \geq 1$ be integers and let $\beta, \xi > 0$ satisfy $\beta
\leq \xi^2 / (3026r^3)$. Let $R$ be a graph over the vertex set $[k]
\times [r]$ such that $\delta(R) > (r-1)k$ and $K_k^r \subset C_k^r
\subset R$. Let $H$ be a graph on $n$ vertices
with maximum degree $\Delta$, and assume that \\
\begin{tabular}{cl}
$(i)$ & $H$ has a labeling of bandwidth at most $\beta n$ and has chromatic number at most $r$. \\
$(ii)$ & For every interval $[a, a + \beta^2 n] \subset [1, n]$, there
exists a vertex $v \in [a, a+\beta^2n]$ such that \\
 & $N_H(v)$ is an independent set. \\
$(iii)$ & $(m_{i,j})_{1 \leq i \leq k, 1 \leq j \leq r}$ is an
$r$-equitable integer partition of $n$ with $m_{i,j} \geq 200\beta
n$ for every \\
& $1 \leq i \leq k$ and $1 \leq j \leq r$.
\end{tabular}

Then there exists a mapping $f: V(H) \rightarrow [k] \times [r]$ and
a set of special vertices $X \subset V(H)$ with the following
properties. \\
\begin{tabular}{cl}
$(a)$ & $|X| \leq kr \xi n$, \\
$(b)$ & the sets $W_{i,j} := f^{-1}(i,j)$ have size $m_{i,j} - \xi n \leq |W_{i,j}| \leq m_{i,j} + \xi n$ for every $i$ and $j$, \\
$(c)$ & for every edge $\{u, v\} \in E(H)$ we have $\{f(u), f(v)\} \in E(R)$,\\
$(d)$ & if $\{u, v\} \in E(H)$ and, moreover, $u$ and $v$ are both in $V(H) \backslash X$, then
$\{ f(u), f(v) \} \in E(K_k^r)$,\\
$(e)$ & $\forall 1 \le i \le k$, $\exists$ at least $\beta^{-1}$ vertices $w \in \big(\cup_{1 \leq j \leq r} W_{i,j}\big) \setminus \big( \cup_{l=0}^{3} N_H^{(l)}(X) \big)$ whose neighborhood \\
& $N_H(w)$ forms an independent set.
\end{tabular}
\end{lemma}
\begin{pf}
The process of finding a map $f$ which satisfies $(a),(b),(c)$,
and $(d)$ can be found in the proof of Lemma 8 in \cite{MR2448444}. We claim that
$(e)$ is also a byproduct of their proof. It suffices to verify that
for all $1 \le i \le k$, there exists an interval of length
at least $2\beta n$ in the set $\cup_{j=1}^r W_{i,j} \setminus \big(\cup_{s=0}^{3} N_H^{(s)}(X) \big)$,
since by condition $(ii)$ this will give at least $\beta^{-1}$ vertices
in this set which have independent neighborhoods. The stronger lower
bound of $m_{i,j} \geq 200\beta n$ that we imposed on top of the conditions of
Lemma 8 in \cite{MR2448444} guarantees
that such an interval always exists. We omit the details.
\end{pf}

Let $G'$ be a given graph and use Lemma \ref{lemma_lemmaforG} to get
a set $B$, a `temporary' partition of $V \backslash B$ which we can adjust (see the
discussion following Lemma \ref{lemma_lemmaforG}), and an integer
partition $(m_{i,j})_{1 \le i \le k, 1 \le j \le r}$. To simplify the
explanation, assume for a moment that the set $B$ is empty.
Use this integer partition $(m_{i,j})_{i,j}$ as an input to Lemma \ref{lemma_lemmaforH},
and we get a partition $(W_{i,j})_{i,j}$ of
the vertex set of $H$, such that the integer partition
$(|W_{i,j}|)_{i,j}$ is close to
$(m_{i,j})_{i,j}$. Thus by Lemma
\ref{lemma_lemmaforG}, we can get a partition $(V_{i,j})$ of $V(G)$
such that $|V_{i,j}| = |W_{i,j}|$ for all $i,j$.

Ideally, we want all the pairs $(V_{i,j}, V_{i', j'})$ to be
super-regular. But in reality, the super-regular pairs are only
guaranteed over $K_k^r$, and the set $X$ in Lemma
\ref{lemma_lemmaforH} is designed to overcome this difficulty.
Observe that all the edges of $H$ which are not incident to $X$
corresponds to $K_k^r$ in the homomorphic image (property (d) of
Lemma \ref{lemma_lemmaforH}). Thus if we can find an embedding of
vertices of $X$ first, so that its neighborhood $Y := N(X)$ is only
`mildly' restricted, then we can extend this embedding by using the
version of the blow-up lemma as in Theorem \ref{thm_blowuplemma}.
The next lemma, which is Lemma 9 in \cite{MR2448444}, can be used to
embed $X$ so that the number of the possible images of each vertex
$y \in Y$ is still large enough.

\begin{lemma} \label{lemma_partialembedding}
For every integer $\Delta \geq 2$ and every $d \in (0,1]$ there
exist constants $c=c(\Delta,d)$ and $\varepsilon_0=\varepsilon_0(\Delta,d)$ such that for every positive
$\varepsilon \leq \varepsilon_0$ the following is true.

Let $R$ be a graph over the vertex set $V(R_k) = [k] \times [r]$ and
$G$ be a graph on $n$ vertices with $V(G) = \bigcup_{1 \leq i \leq
k, 1 \leq j \leq r} V_{i,j}$, such that $|V_{i,j}| \geq (1 -
\varepsilon)n/(kr)$ for all $1 \leq i \leq k, 1 \leq j \leq r$ and
as a partition, $(V_{i,j})$ is $(d,\varepsilon)$-regular on $R$.
Furthermore, let $\Gamma$ be a graph with $V(\Gamma) = X \cup Y$ and
$f: V(\Gamma) \rightarrow V(R)=[k] \times [r]$ be a mapping with
$\{f(a), f(a')\} \in E(R)$ for all $\{a,a'\} \in E(\Gamma)$.

If $|V(\Gamma)| \leq \varepsilon_0 n / (kr)$ and $\Delta(\Gamma)
\leq \Delta$, then there exists an injective mapping $g: X
\rightarrow V(G)$ with $g(x) \in V_{f(x)}$ for all $x \in X$ such
that for all $y \in Y$ there exist sets $C_y \subset V_{f(y)}
\backslash g(X)$ such that \\
\begin{tabular}{cl}
$(i)$ & $g$ is a graph homomorphism of $\Gamma[X]$ to $G$,\\
$(ii)$ & for all $y \in Y$ we have $C_y \subset N_G(g(x))$ for all $x \in N_\Gamma(y) \cap X$, and \\
$(iii)$ & $|C_y| \geq c|V_{f(y)}|$ for every $y \in Y$.
\end{tabular}
\end{lemma}

\section{Technical lemmas}
\label{section_thekeylemma}

In this section we prove Lemma \ref{lemma_lemmaforG} by using the
following useful statement. This statement hints where the set $B$ in Lemma
\ref{lemma_lemmaforG} comes from.
\begin{lemma}\label{lemma_badset}
Let $0 < p \leq 1$ be fixed and $T$ be an integer. Then for every
$\varepsilon>0$, there exists a constant $b_0=b_0(p,T,\varepsilon)$
such that $G=G(n,p)$ a.a.s.~satisfies the following. For arbitrary
subsets $V_1, \ldots, V_T$ of the vertex set $V$ with $|V_i| \geq
\varepsilon n$ for all $1 \le i \le T$, there exists a set $B$ of
size at most $b_0$ such that for all $v \in V \setminus B$, we have
$\deg(v, V_i) \in [(1 - \varepsilon) |V_i|p, (1 +
\varepsilon)|V_i|p]$ for all $1 \le i \le T$.
\end{lemma}
\begin{pf}
Let $b'$ be a constant to be chosen later. As a first step, we fix a
set $W \subset V$ of size at least $\varepsilon n$, and analyze the
probability of there being $b'$ vertices $v$ such that $\deg(v, W)
\notin [(1 - \varepsilon) |W|p, (1 + \varepsilon)|W|p]$. Let $B$ be
a set of size $b'$ and assume that for all $v \in B$, we have
$\deg(v, W) < (1 - \varepsilon) |W|p$. Then by definition, $e(B, W)
< |B| \cdot (1 - \varepsilon)|W|p$. We estimate the probability of
this event. Note that $B$ is a set of constant size and $W$ has size
$|W| \ge \varepsilon n$, thus it suffices to bound the probability
of $e(B, W \setminus B) < |B| \cdot (1 - \varepsilon/2)|W \setminus
B|p$. Since $e(B, W \setminus B)$ has expectation $|B||W\setminus
B|p$ and is a sum of independent binomial random variables, we can
use Chernoff inequality to get,
\[ P\big( e(B, W \setminus B) < (1 - \varepsilon/2)|B||W\setminus B|p \big)
\le e^{-\Omega_{\varepsilon}(b'np)}.
\]
Thus for a fixed set $B$ of size $b'$ and $W$ of size at least
$\varepsilon n$, the probability that all the vertices $v \in B$
have $\deg(v, W) < (1 - \varepsilon) |W|p$ is
$e^{-\Omega_{\varepsilon}(b'np)}$. Take the union bound of this
event over all choices of $B$ and $W$ and we can conclude that the
probability of there existing such sets $B$ and $W$ in $G$ is at
most ${n \choose b'} \cdot 2^n \cdot e^{-\Omega_{\varepsilon}(b'np)}
= o(1)$ as long as $b'=b'(\varepsilon,p)$ is large enough. In other
words, a.a.s.~every set $W$ of size $|W| \ge \varepsilon n$ has at most
$b'$ vertices $v$ such that $\deg(v, W) < (1 - \varepsilon) |W|$.

Given subsets $V_1, \ldots, V_T$ of size at least $\varepsilon n$,
the previous observation implies that there are at most $b'T$
vertices which have $\deg(v, V_i) < (1 - \varepsilon) |V_i|$ for
some $1 \le i \le T$, and similarly at most $b'T$ vertices which
have $\deg(v, V_i) > (1 + \varepsilon) |V_i|$ for some $1 \le i \le
T$. Therefore by setting $b_0 = 2b'T$, we can derive the conclusion
of the lemma.
\end{pf}

The proof of Lemma \ref{lemma_lemmaforG} consists of two steps. The
first step is to show the existence of a `temporary' partition
$(U_{i,j})_{1 \leq i \leq k, 1 \leq j \leq r}$ which has size
$m_{i,j} := |U_{i,j}|$ for all $i,j$ (see the discussion following
the statement of Lemma \ref{lemma_lemmaforG}). Once this partition
is constructed, we select sets $V_{i,j}^*$ arbitrarily within the
`temporary' set $U_{i,j}$, and for a given integer partition
$(n_{i,j})$, modify the partition slightly without moving the
vertices in $V_{i,j}^*$ to make the sizes of the partition as
desired.

The two lemmas below establish stability results for regular and
super regular pairs. They basically say that  regularity can be changed into
super-regularity by small perturbation (Lemma
\ref{lemma_preserveregularity}), and regularity and
super-regularity are stable under small perturbation (Lemma
\ref{lemma_preservesuperregularity}). These can be found in
Proposition 13, and 14 of \cite{MR2448444}.

\begin{lemma} \label{lemma_preserveregularity}
Fix $\varepsilon,d>0$. For any graph $G$ and $\varepsilon$-regular
partition $V_1, \ldots, V_k$ with $(d,\varepsilon)$-reduced graph
$R$, let $S$ be a subgraph of $R$ with $\Delta(S) \leq \Delta$. Then
for each vertex $i$ of $S$, we can find a set $V_i' \subset V_i$ of
size $(1- \varepsilon\Delta)|V_i|$ such that for every edge $\{i,
j\} \in E(S)$ the pair $(V_i', V_j')$ is $(d -
\varepsilon(\Delta+1), \varepsilon/(1-\varepsilon
\Delta))$-super-regular. Moreover, for every edge $\{i, j\}$ of the
original reduced graph $R$, the pair $(V_i', V_j')$ is still $(d -
\varepsilon(\Delta+1), \varepsilon/(1-\varepsilon \Delta))$-regular.
\end{lemma}
\begin{lemma} \label{lemma_preservesuperregularity}
Let $(A,B)$ be an $(d,\varepsilon)$-regular pair and let $(\hat{A},
\hat{B})$ be a pair such that $|\hat{A} \Delta A | \leq
\hat{\alpha}|\hat{A}|$ and $|\hat{B} \Delta B | \leq
\hat{\beta}|\hat{B}|$ for some $0 \leq \hat{\alpha}, \hat{\beta}
\leq 1$. Then, $(\hat{A}, \hat{B})$ is an $(\hat{d},
\hat{\varepsilon} )$-regular pair with $\hat{d} := d -
2(\hat{\alpha} + \hat{\beta})$ and $\hat{\varepsilon} := \varepsilon
+ 3(\hat{\alpha}^{1/2} + \hat{\beta}^{1/2})$. If, moreover, $(A,B)$
is $(d,\varepsilon)$-super-regular and each vertex $v$ in $\hat{A}$
has at least $d|\hat{B}|$ neighbors in $\hat{B}$ and each vertex $v$
in $\hat{B}$ has at least $d|\hat{A}|$ neighbors in $\hat{A}$, then
$(\hat{A}, \hat{B})$ is $(\hat{d},\hat{\varepsilon})$-super-regular.
\end{lemma}

Next lemma is an immediate corollary of the bandwidth theorem proved in \cite{MR2448444}
(which is Theorem 1 there).

\begin{lemma} \label{lemma_backbonelemma}
Given an integer $r \geq 1$ and a constant $\gamma > 0$,
any sufficiently large graph $G$
on $n$ vertices with minimum degree $(1 - 1/r + \gamma)n$
contains a copy of $C_m^r$ with $m = \floor{n/r}$.
\end{lemma}

Now we are ready to prove Lemma \ref{lemma_lemmaforG}.\\

\begin{pfof}{Lemma \ref{lemma_lemmaforG}}
Given $r \ge 2, p, \gamma$, choose $d \leq \gamma p /90$ and let
$\varepsilon_0 = \min \{
\varepsilon_{\ref{lemma_reducedgraphmindegree}}(p, r, \gamma),
d/(12r) \}$. Assume that an $\varepsilon \leq \varepsilon_0$ is
given, and let $\varepsilon' = \gamma p\varepsilon^6 / (1152r)$
and $d' = d + 2\varepsilon$. Let $t = \max\{ 4r / \gamma,
1/(\varepsilon')\}$ and $T = T_{\ref{thm_regularitylemma}}(t,
\varepsilon')$. Let $b_0 = b_{\ref{lemma_badset}}(p, T,
\varepsilon')$.

Let $G=G(n,p)$ and let $G' \subset G$ be a subgraph with $\delta(G')
\geq (1-1/r + \gamma)np$. By using the degree form of the regularity
lemma (Lemma \ref{thm_regularitylemma}), we obtain a graph $G''
\subset G'$ and a pure $(d',\varepsilon')$-regular partition
$(U_i)_{0 \leq i \leq s}$ of $G''$ with reduced graph $R$ and $t
\leq s \leq T$. From now on we will only consider the graph $G''$,
unless mentioned otherwise. Remove at most $r-1$ parts and put them
into the set $U_0$ so that we can assume $s = kr$ for some integer
$k$. Note that by Lemma \ref{thm_regularitylemma} $(ii)$,
\[ \delta(G'') \geq \delta(G') - (d' + \varepsilon')n \ge (1-1/r + \gamma)np - (d'+\varepsilon')n \ge (1-1/r + 4\gamma/5)np, \]
and thus by Lemma \ref{lemma_reducedgraphmindegree} we have $\delta(R) \geq (1 - 1/r +
\gamma/2)s$. Let $m := |U_i|$ and note $|U_0| \le \varepsilon' n + (r-1)\frac{n}{t} \le r \varepsilon' n$,
so we have $ms = mkr \le n \le mkr/(1-r\varepsilon')$.

By Lemma \ref{lemma_backbonelemma}, $R$ contains a copy of $C_k^r$.
Thus we may assume that $R$ is a graph over the vertex set $[k]
\times [r]$ with $K_k^r \subset C_k^r \subset R$. Rename the parts
$U_x$ as $U_{i,j}$ according to this new vertex set of $R$ to get a
vertex partition $U_0 \cup \bigcup_{1 \leq i \leq k, 1 \leq j\leq r}
U_{i,j}$. Then by applying Lemma \ref{lemma_preserveregularity} with
$S = K_k^r$ and $\Delta = r-1$, one can obtain a new partition $U_0'
\cup \bigcup_{1 \leq i \leq k, 1 \leq j\leq r} U_{i,j}'$ which is
$(d' - \varepsilon' r,
\varepsilon'/(1-\varepsilon'r))$-super-regular on $K_k^r$, $(d' -
\varepsilon'r, \varepsilon'/(1 - \varepsilon'r))$-regular on $R$,
and $|U_{i,j}'| = (1 - \varepsilon'r)m$. Since all the discarded
vertices of $U_{i,j}$ are collected into $U_0'$, we have
 $|U_0'| \leq |U_0| + \varepsilon'mkr^2 \le \varepsilon' r n + \varepsilon'mkr^2 \leq 2\varepsilon'r n$. Applying Lemma
\ref{lemma_badset} to the sets $U_{i,j}'$, we get a set $B$ such that
for all $v \in V \backslash B$, $\deg_{G''}(v, U_{i,j}') \le \deg_{G}(v, U_{i,j}') \leq (1
+\varepsilon') mp$ for all $i \in [k], j \in [r]$. Remove all the
vertices of $B$ belonging to $U_{i,j}'$ for $i \in [k], j \in [r]$, and put it into
$U_0'$, and then remove some more vertices from each partition so
that the number of vertices in each part is the same for all $i,j$. Since $B$ is a set of
constant size, asymptotically the effect of this process is
negligible and we may use the same bounds on the size of the sets as
before.

We would like to spread the vertices in the exceptional set $U_0'
\backslash B$ into $(U_{i,j}')_{1 \leq i \leq k, 1 \leq j \leq r}$
while keeping the $r$-equitable property of the partition,
regularity on $R$ and super-regularity on $K_k^r$. For a vertex $u
\in U_0' \setminus B$ call an index $i$ good if $u$ has at least $d'm$ neighbors
in each $U_{i,j}'$ for all $j \in [r]$. Let $g_u$ be the number of good
indices for $u$, and let $U_i' = \bigcup_{1 \leq j \leq r}
U_{i,j}'$. Note that if $i$ is a good index for $u$, then we can add
$u$ to any part of $U_i'$ without destroying the super-regularity on
$K_k^r$. By the definition of $B$, for $u \in V \backslash B$ and
arbitrary $i \in [k], j \in [r]$, $\deg_{G''} (u, U_{i,j}') \le (1 +
\varepsilon')mp$ and so $\deg_{G''} (u, U_i') \leq (1 +
\varepsilon')rmp$ in general. However, if $i$ is not a good index
for $u$, then $u$ can only have at most $d'm$ neighbors in one of
the parts, and we have the bound $\deg_{G''} (u, U_i') \leq
(1+\varepsilon')(r-1)mp + d'm$. Thus we have
\[ \deg_{G''}(u, U_1' \cup \ldots \cup U_k') \le g_u (1+\varepsilon')rmp + (k - g_u) \big( (1+\varepsilon')(r-1)mp + d'm \big). \]
On the other hand, since $G''$ has minimum
degree at least $(1-1/r + 4\gamma/5)np$, and $|U_0'| \leq
2\varepsilon' r n$, we have,
\[ \deg_{G''}(u, U_1' \cup \ldots \cup U_k') = \deg_{G''}(u, V \backslash U_0') \ge (1 - 1/r + 4\gamma/5)np - 2\varepsilon' rn \ge (1 - 1/r + 3\gamma/4)np. \]
Combine these bounds to get,
\[ (1 - 1/r + 3\gamma/4)np \le g_u (1+\varepsilon')rmp + (k - g_u) \big( (1+\varepsilon')(r-1)mp + d'm \big). \]
Using the fact $mkr \le n$, we can divide the left hand side by $np$, and
right hand side by $mkrp$ to get,
\begin{align*}
1 - 1/r + 3\gamma/4 &\leq \frac{g_u}{k} (1 + \varepsilon') + \left(1 - \frac{g_u}{k}\right)\left(\left(1 - \frac{1}{r}\right)(1 + \varepsilon') + \frac{d'}{rp} \right) \\
&\le \frac{g_u}{kr} (1 + \varepsilon') + \left(1 - \frac{1}{r}\right) + \varepsilon' + \frac{d'}{rp},
\end{align*}
which implies $g_u \geq \gamma kr /2$. Pick a vertex in $U_0'
\setminus B$ one by one, and assign one of its good index to it as
follows. Always pick the index which has been assigned the least
number of vertices so far. In this way, we can assign an index to
every vertex $U_0' \setminus B$ so that each index gets assigned at
most $2|U_0'|/(\gamma kr)$ vertices. By using the fact $|U_0'| \le 2
\varepsilon'r n$ and $n \le mkr / (1 -  r\varepsilon')$ we get,
\[ \frac{2|U_0'|}{\gamma kr} \leq \frac{4\varepsilon'rn}{\gamma kr} \leq \frac{4\varepsilon'mkr^2}{(1 - r\varepsilon')\gamma kr} \le \frac{8\varepsilon'r}{\gamma}m \le \frac{\varepsilon^6}{144}m =: \alpha m.\]
For each index $i$, spread the vertices of $U_0'$ assigned to it as
evenly as possible into $U_{i,j}'$ for $j \in [r]$ so that the
resulting partition $(U_{i,j}'')_{1 \le i \le k, 1 \le j \le r}$ is
$r$-equitable. Recall that (i) all the vertices assigned to an index
have degrees at least $d'm$ in every part belonging to that index,
and (ii) $(U_{i,j}')_{1 \le i \le k, 1 \le j \le r}$ was $(d' -
\varepsilon'r, \varepsilon/(1-\varepsilon'r))$-super-regular on
$K_k^r$ and $(d' - \varepsilon'r,
\varepsilon/(1-\varepsilon'r))$-regular on $R$. Furthermore, the
sets $U_{i,j}'$ had size $(1 - \varepsilon'r)m$, and $|U_{i,j}''
\Delta U_{i,j}'| \le \left\lceil \alpha m /r \right\rceil \le \alpha
m \le \alpha|U_{i,j}'|/(1-\varepsilon'r) \le 2 \alpha |U_{i,j}'|$.
Thus by Lemma \ref{lemma_preservesuperregularity} we know that
$(U_{i,j}'')_{1 \le i \le k, 1 \le j \le r}$ is $(d' - \varepsilon'
r - 8\alpha, \varepsilon'/(1-\varepsilon'r) + 6\sqrt{2}
\alpha^{1/2})$-super-regular on $K_k^r$ and $(d' - \varepsilon' r -
8\alpha, \varepsilon'/(1 - \varepsilon'r) + 6\sqrt{2}
\alpha^{1/2})$-regular on $R$. By the choice of the parameters, we
have,
\[ d' - \varepsilon'r - 8\alpha \geq d + 2\varepsilon - \frac{\gamma p \varepsilon^6}{1152} - \frac{\varepsilon^6}{18} \ge d + \varepsilon, \qquad \textrm{and} \qquad
\frac{\varepsilon'}{1-\varepsilon'r} + 6\sqrt{2} \alpha^{1/2} \leq
2\varepsilon' + \frac{\varepsilon^3}{\sqrt{2}} \le \varepsilon^3. \]
Therefore $(U_{i,j}'')$ is
$(d+\varepsilon,\varepsilon^3)$-super-regular on $K_k^r$ and
$(d+\varepsilon,\varepsilon^3)$-regular on $R$.

Let $m_{i,j} := |U_{i,j}''|$ and note that this satisfies \[ m_{i,j}
\geq |U_{i,j}'| \ge (1 - \varepsilon' r)m \geq \frac{(1 -
\varepsilon' r)^2n}{kr} \ge \frac{(1 - \varepsilon)n}{kr} \] for all
$i \in [k], j \in [r]$. Then fix an arbitrary set $V_{i,j}^* \subset
U_{i,j}''$ of size $(1-3\varepsilon^3 r)m_{i,j}$ for all $i \in [k],
j \in [r]$ and note that $(1 - 3\varepsilon^3 r)m_{i,j} \ge (1 -
\varepsilon)m_{i,j}$ so that $(iv)$ holds. Since $|U_{i,j}'' \Delta
V_{i,j}^*| = 3\varepsilon^3 r |U_{i,j}''|$, by Lemma 4.3, the
partition $(V_{i,j}^*)$ will be $(d+\varepsilon - 12\varepsilon^3 r,
\varepsilon^3 + 6\sqrt{3}(\varepsilon^3 r)^{1/2})$-regular on $R$,
and in particular $(d, \varepsilon)$-regular on $R$. This concludes
the first part of Lemma \ref{lemma_lemmaforG} where given a graph
$G'$, we obtain a subgraph $G''$, a set $B$, sets $(V_{i,j}^*)$ which are
$(d,\varepsilon)$-regular on $R$, and an $r$-equitable integer
partition $(m_{i,j})$ of $n - |B|$.

It remains to show that given another integer partition $(n_{i,j})$,
we can find a partition $(V_{i,j})$ of $V \setminus B$ with
$|V_{i,j}| \ge n_{i,j}$ for all $i,j$. This partition will be
obtained from the partition $(U_{i,j}'')$ by {\em pushing around the
vertices}. This is a process of moving vertices from one partition
to another while keeping regularity and super-regularity of pairs.
For example, say that we want to move one vertex from $U_{1,1}''$ to
$U_{2,1}''$. Then by the regularity of $(U_{i,j}'')$ on $C_{k}^r$,
there exists a vertex $u \in U_{1,1}''$ which has high degree in all
the sets $U_{2,j}''$ for $2 \le j \le r$. Moving this vertex to
$U_{2,1}''$ will not destroy the regularity and super-regularity of
pairs. One must observe that the proof in \cite{MR2448444} allows to
fix a set $V_{i,j}^*$ of size $(1-3\varepsilon^3 r)m_{i,j}$ and
always choose a vertex outside of it to push around (this follows
from the $(d+\varepsilon,\varepsilon^3)$ regularity on $R$). After
the process of pushing around the vertices is done, the size of the
sets $U_{i,j}''$ will change, and thus affect the super-regularity
and regularity between parts. This is where we want to choose $\xi_0
= \xi_0(r,p,\gamma,\varepsilon)$ to be small enough. By doing so, we
can make sure that the sets $U_{i,j}''$ changes only by some small
amount, and since $(U_{i,j}'')$ is
$(d+\varepsilon,\varepsilon^3)$-super-regular on $K_k^r$ and
$(d+\varepsilon,\varepsilon^3)$-regular on $R$, the resulting
partition $(V_{i,j})$ will still be $(d, \varepsilon)$-super-regular
and regular on $K_k^r$ and $R$, respectively. For further details,
we refer the reader to the proof of Lemma 6 in \cite{MR2448444}.
\end{pfof}

\section{Main Theorem}
\label{section_maintheorem}

In this section we prove the main theorem.

\begin{thm} \label{thm_mainthm}
For fixed integers $r, \Delta$, and reals $0 < p \leq 1$ and $\gamma
> 0$, there exists a constant $\beta > 0$ such that a.a.s., any
spanning subgraph $G'$ of $G(n,p)$ with minimum degree $\delta(G')
\ge (1 - 1/r + \gamma)np$ contains every $n$-vertex graph $H$ which
satisfies the following properties. $(i)$ $H$ is $r$-chromatic,
$(ii)$ has maximum degree at most $\Delta$, $(iii)$ has bandwidth at
most $\beta n$ with respect to a labeling of vertices by
$1,2,\ldots, n$, and $(iv)$ for every interval $[a, a+\beta^2 n]
\subset [1,n]$, there exists a vertex $v \in H$ such that $N_H(v)$
is an independent set.
\end{thm}
\begin{pf}
First we will adjust the parameters. We may assume that $r \ge 2$,
since the case $r=1$ is trivial. Given $r,\Delta, p, \gamma$, take
$d=d_{\ref{lemma_lemmaforG}}(r,p,\gamma)$, $c =
\min\{c_{\ref{lemma_partialembedding}}(\Delta, d/2), (d/8)^\Delta
\}$, and $\alpha = \alpha_{\ref{thm_blowuplemma}}(d/2,\Delta,c,r)$.
Then let
\[ \varepsilon = \frac{1}{2} \min \left\{\varepsilon_{\ref{thm_blowuplemma}}(\frac{d}{2},\Delta,c,r), \varepsilon_{\ref{lemma_partialembedding}}(\Delta,\frac{d}{2}), \varepsilon_{\ref{lemma_lemmaforG}}(r,p,\gamma), \frac{dp}{6r\Delta}, \left(\frac{d}{8}\right)^\Delta \right\}, \]
$b_0 = b_{\ref{lemma_lemmaforG}}(r,p, \gamma, \varepsilon)$, $K_0 = K_{\ref{lemma_lemmaforG}}(r,p,\gamma,\varepsilon)$, and
\[ \xi = \frac{1}{2}\min \left\{\xi_{\ref{lemma_lemmaforG}}(r,p,\gamma,\varepsilon), \frac{(1-\varepsilon) \alpha \varepsilon^2 c}{144\Delta (K_0r)^2} \right\}. \]
Finally, choose $\beta \leq \min \{ \xi^2/(6052r^3), 1/(b_0 \Delta^5)\}$.

Lemma \ref{lemma_lemmaforG} applied to $G'$ provides us a subgraph
$G'' \subset G'$, a graph $R$ over the vertex set $[k] \times [r]$
with $k \le K_0$, a set $B$ with $|B| = b \leq b_0$, sets $(
V_{i,j}^* )_{1 \leq i \leq k, 1 \leq j \leq r}$, and a $r$-equitable
integer partition $(m_{i,j})_{1 \leq i \leq k, 1 \leq j \leq r}$
satisfying $(i), (ii), (iii), (iv)$. Given this partition
$(m_{i,j})$, apply Lemma \ref{lemma_lemmaforH} to $H$ and get a
partition $W_{i,j}$ of $H$ satisfying (a),(b),(c),(d),(e) of the
lemma. Since an embedding of $H$ into $G''$ is also an embedding
into $G'$, by abusing notation, we will denote $G'$ for the graph
$G''$. Note that by doing this, we can only guarantee $\delta(G')
\ge (1 - 1/r + 4\gamma/5)np$.

To control the set $B$, we will find vertices of $H$ which can be
mapped into the set $B$. Note that for this step, the set $B$
contained in $V(G')$ comes first, and then we look at $H$ to decide
which of its vertices can be mapped into $B$. Considering the fact
that we are trying to embed a particular given graph $H$ into $G'$,
this step might seem somewhat peculiar.
\begin{claim} \label{clm_verticesforB}
There exists a set $Z \subset V(H) \setminus \left(\cup_{s=0}^{2}
N^{(s)}(X)\right)$, and a one-to-one graph homomorphism $g : Z
\rightarrow V(G')$ which satisfies the following properties.
\begin{enumerate}[(i)]
  \setlength{\itemsep}{1pt}
  \setlength{\parskip}{0pt}
  \setlength{\parsep}{0pt}
\item $B \subset g(Z) \subset B \cup (\cup_{i,j}V_{i,j}^*)$,
\item for $W_B = g^{-1}(B)$, $Z = W_B \cup N_H(W_B)$.
\item for $w \in N_H^{(2)}(W_B)$, assume that $w \in W_{i,j}$.
Then there exists a set $C_w \subset V_{i,j}^* \setminus g(Z)$
of size $|C_w| \ge 2c m_{i,j}$ which is contained in the common
neighborhood of all vertices in $g(N_H(w) \cap Z)$.
\end{enumerate}
\end{claim}
The proof of this claim will be given later. Once we apply this
claim, we obtain a partial embedding $g$ of $H$ which embeds the
vertices $Z$, and constrains the image of every vertex $w \in N_H(Z)
\setminus Z$ to some set $C_w$. Moreover, the set $B$ is covered by
the image of this map.

Next, we adjust the partition of $G'$ in order to embed the
remaining vertices of $H$. The goal is to obtain a partition
in which the sets $V_{i,j}$ have size $n_{i,j} = |W_{i,j} \backslash Z| +
|V_{i,j}^* \cap g(Z)|$, where the first term comes from the number
of remaining vertices to be mapped and the second term comes from
the vertices which have already been mapped to $V_{i,j}^*$. Let $\delta_{i,j} =
|V_{i,j}^* \cap g(Z)|$, and note that $\sum_{1 \leq i \leq k, 1 \leq j \leq r}
\delta_{i,j} \leq |g(Z)| = |Z| \le (\Delta+1)b_0$ by part $(ii)$ of Claim \ref{clm_verticesforB}.
Since $\Delta, b_0$ are constants and $m_{i,j}$ is linear in $n$ for all $i,j$,
\begin{align*}
 n_{i,j} &\le |W_{i,j}| + \delta_{i,j} \le (1 + \xi)m_{i,j} + (\Delta+1)b_0 \le (1 + 2\xi)m_{i,j} , \quad \textrm{and} \\
 n_{i,j} &\ge |W_{i,j}| - |Z| \ge |W_{i,j}| - (\Delta+1)b_0 \geq (1-\xi)m_{i,j} - (\Delta+1)b_0 \ge (1 - 2\xi)m_{i,j}.
\end{align*}
Therefore $n_{i,j} \in [(1 - \xi_{\ref{lemma_lemmaforG}})m_{i,j}, (1
+ \xi_{\ref{lemma_lemmaforG}})m_{i,j}]$. Moreover, we have
\[ \sum_{i,j}n_{i,j} = \sum_{i,j}|W_{i,j} \backslash Z| +
|V_{i,j}^* \cap g(Z)| = \left(\sum_{i,j}|W_{i,j}|\right) - |Z| + (|g(Z)| - |B|) = n - b. \]

Thus we can use Lemma \ref{lemma_lemmaforG} to obtain a partition
$(V_{i,j})_{1 \leq i \leq k, 1 \leq j \leq r}$ of the vertices $V
\backslash B$ such that $|V_{i,j}| = n_{i,j}$ for all $i,j$,
$(V_{i,j})$ is $(d,\varepsilon)$-regular on $R$, and
$(d,\varepsilon)$-super-regular on $K_k^r$. Then since $g(Z) \subset
V_{i,j}^* \subset V_{i,j}$, by defining $V_{i,j}' = V_{i,j}
\backslash g(Z)$, we have, $|V_{i,j}'| = n_{i,j} - \delta_{i,j} =
|W_{i,j} \backslash Z|$. Note that we removed only at most constant
number of vertices from $V_{i,j}$ to obtain $V_{i,j}'$. Thus by
Lemma \ref{lemma_preservesuperregularity}, $(V_{i,j}')_{1 \leq i
\leq k, 1 \leq j \leq r}$ is $(d-\varepsilon,2\varepsilon)$-regular
on $R$ and $(d-\varepsilon,2\varepsilon)$-super-regular on $K_k^r$.
Let $V' := \bigcup_{1\le i \le k, 1 \le j \le r} V_{i,j}'$. Since $d
- \varepsilon \ge d/2$, we may assume that the partition
$(V_{i,j}')$ is $(d/2, 2\varepsilon)$-regular and super-regular,
respectively.

We would like to find an embedding of the remaining vertices of $H$
so that $W_{i,j} \setminus Z$ gets mapped to $V_{i,j}'$ for all
$i,j$, and every vertex $w \in N(Z) \setminus Z$ gets mapped to a
vertex in $C_w$. Recall that $X$ is a subset of $V(H)$ obtained in
Lemma \ref{lemma_lemmaforH} and $|X \cup N(X)| \le (\Delta + 1)kr\xi n \le (\varepsilon_{\ref{lemma_partialembedding}}/(kr))n$. Apply Lemma
\ref{lemma_partialembedding} with the set $X$ and $Y = N(X)
\setminus X$ to embed the vertices $X$ into $V'$ so that $(i)$,
$(ii)$, $(iii)$ of Lemma \ref{lemma_partialembedding} holds. Now we
have a new set of constraints, namely, every $y \in Y$ has a set
$C_y$ which it has to be mapped to. Since  $Z \subset V(H) \setminus
\left(\cup_{s=0}^{2} N^{(s)}(X)\right)$, the set $Y$ and $N(Z)
\setminus Z$ are disjoint, thus the constraints coming from the
vertices $Y$ and the ones coming from $N(Z) \setminus Z$ will not
interfere with each other. Extend the map $g$ which embedded the
vertices $Z$ so that $g$ is an embedding of $X \cup Z$. Let
$V_{i,j}'' := V_{i,j}' \setminus g(X) = V_{i,j} \setminus g(X \cup
Z)$ and $V'' = \bigcup_{i,j}V_{i,j}''$. Then by $m_{i,j} \ge (1 - \varepsilon)n / (kr)$ from Lemma
\ref{lemma_lemmaforG} and $|V_{i,j}| = n_{i,j} \ge (1 -
2\xi)m_{i,j}$,
\[ |V_{i,j} \setminus V_{i,j}''| \le |X| + |Z| \le kr\xi n + (\Delta + 1)b_0 \le \frac{2\xi k^2r^2 m_{i,j}}{1 - \varepsilon} \le \frac{2\xi k^2r^2}{(1 - 2\xi)(1 - \varepsilon)} |V_{i,j}| \le \frac{\varepsilon^2}{36} |V_{i,j}|. \]
Recall that the partition $(V_{i,j})$ was $(d, \varepsilon)$-regular
on $R$ and $(d,\varepsilon)$-super-regular on $K_k^r$. Consequently,
by Lemma \ref{lemma_preservesuperregularity} with $\hat{\alpha} =
\hat{\beta} = \varepsilon^2/36$, the partition $(V_{i,j}'')$ is
$(d-\varepsilon^2/9,\varepsilon + \varepsilon)$-regular on $R$ and
 $(d-\varepsilon^2/9,\varepsilon +
\varepsilon)$-super-regular on $K_k^r$. We may assume that
$(V_{i,j}'')$ is $(d/2, 2\varepsilon)$-regular and super-regular,
respectively.

Let $f$ be the graph homomorphism of $H$ to $R$ given in Lemma
\ref{lemma_lemmaforH}. Since we finished embedding $X$, by (d) of
Lemma \ref{lemma_lemmaforH}, the homomorphic image under $f$ of all
the remaining edges of $H$ correspond to $K_k^r$ in the graph $R$. Thus
once we check that the parameters are chosen correctly, we can apply
the blow-up lemma, Theorem \ref{thm_blowuplemma}, to each of the
partition $(V_{i,j}'')_{1 \leq j \leq r}$ for fixed $i \in [k]$
separately, to find an embedding of the remaining vertices $V(H)
\backslash (X \cup Z)$ which is consistent with the map $g$.

In the remaining part of the proof, we verify that the parameters
are chosen so that we can apply the blow-up lemma. The previously
embedded vertices constrains the possible images of vertices in
$N_H(Z)\setminus Z$ and $Y = N_H(X) \setminus X$. For a vertex $w
\in N_H(Z) \setminus Z$, by Claim \ref{clm_verticesforB}, the image
of $w$ were constrained to a set $C_w \subset V_{i,j}^*$ of size at
least $2c m_{i,j}$ for some $i,j$. Among these vertices, some could
have been used for the sets $X$, but the number of remaining
vertices in $C_w$ is still at least
\[ 2c m_{i,j} - |X| \ge 2c m_{i,j} - k r \xi n \ge 2c m_{i,j} - \frac{(kr)^2\xi }{(1 - \varepsilon)}m_{i,j} \ge 2c m_{i,j} - \frac{c}{4} m_{i,j} \ge c n_{i,j} \ge c |V_{i,j}''|, \]
where we used $m_{i,j} \ge (1 - \varepsilon)n / (kr)$ from Lemma
\ref{lemma_lemmaforG} $(ii)$, and $n_{i,j} \le (1 + 2\xi)m_{i,j}$
which we established above, and $n_{i,j} = |V_{i,j}| \ge
|V_{i,j}''|$. For a vertex $y \in Y$, the size of the set $C_y$ is
at least $c |V_{i,j}'| \ge c |V_{i,j}''|$ for corresponding $i,j$ by
Lemma \ref{lemma_partialembedding}.

Moreover, by the choice of $\xi$ depending on $\alpha$, we have
$|N_H(X)| \leq \Delta |X| \leq \Delta k r \xi n \leq (\alpha/4)
m_{i,j} \leq (\alpha/2) n_{i,j}$ for arbitrary $i,j$, and so the
size of $Y$ is less than $(\alpha/2) \min_{i,j} n_{i,j}$. Also,
$N(Z) \setminus Z$ has size at most $|N^{(2)}(W_B)| \le b_0
\Delta^2$ which is a constant. Thus there are at most $\alpha
\min_{i,j} n_{i,j}$ vertices vertices inside $V''$ whose images are
constrained. Finally, note that we picked $2\varepsilon \le
\varepsilon_{\ref{thm_blowuplemma}}(d/2,\Delta,c,r)$, so that $(d/2,
2\varepsilon)$-super-regularity over $K_k^r$ suffices for the
application of the blow-up lemma, Theorem \ref{thm_blowuplemma}.
Once we apply the blow-up lemma, we can find a mapping which embeds
all the remaining vertices of $H$, and when combined with the
previous mappings, forms a graph homomorphism of $H$ into $G$.
\end{pf}

\begin{pfof}{Claim \ref{clm_verticesforB}}
For a vertex $v \in B$, since $n = |B| + \sum_{i,j} m_{i,j}$,
\[ \Big| V \backslash \big( \cup_{1 \leq i \leq k, 1 \leq j \leq r} V_{i,j}^* \big) \Big| = n - \sum_{i,j}|V_{i,j}^*| = |B| + \sum_{i,j}(m_{i,j} - |V_{i,j}^*|) \le b_0 + \sum_{i,j} \varepsilon m_{i,j} \le b_0 + \varepsilon n, \]
and $v$ has at least
\[ \delta(G') - (b_0 + \varepsilon n) \ge (1 - 1/r)np  - 2\varepsilon n = (1 - 1/r - 2\varepsilon p^{-1})np \]
neighbors in $\cup_{i,j} V_{i,j}^*$. By the fact $\sum_{i,j}m_{i,j}
\le n$, this implies that there exists an index $(s,t)$ such that
$v$ has at least $(1 - 1/r - 2\varepsilon p^{-1})m_{s,t}p \ge
\frac{1}{3}m_{s,t}p$ neighbors in $V_{s,t}^*$. Since $R$ has $rk$
vertices and $\delta(R)
> (r-1)k$ by Lemma \ref{lemma_lemmaforG} $(i)$, by pigeonhole
principle there exists an index $s' \in [k]$ such that $(s, t)$ is
adjacent to $(s',j)$ in $R$ for all $j \in [r]$. By property $(e)$
of Lemma \ref{lemma_lemmaforH} there exists at least $1/\beta$
vertices in $\big( \cup_{1 \leq j \leq r}W_{s',j} \big) \backslash
\big( \cup_{l=0}^{3} N^{(l)}(X) \big)$ which have independent
neighborhoods. Since $|B|\Delta^5 \leq b_0 \Delta^5 \leq 1/\beta$,
we can assign one such vertex $h_v$ to each $v \in B$ so that the
vertices $h_v$ have distance at least 5 to each other in $H$ (we
want them to be far apart from each other so that later they do not
constrain the same set of vertices). Thus we have assigned $g(h_v) =
v$ (see figure \ref{figure_verticesforB}).
\begin{figure}
\centering
\includegraphics[width=14.0cm]{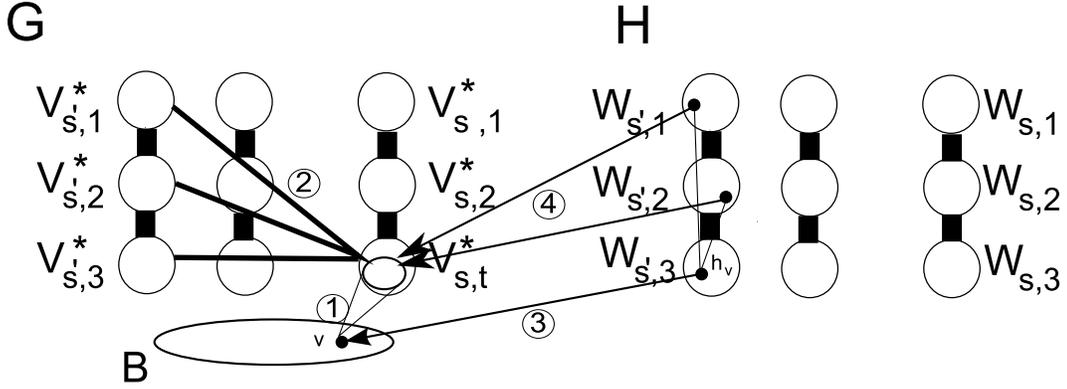}
\caption{Assigning vertices to $B$. Number indicates the logical order.}
\label{figure_verticesforB}
\end{figure}

Fix a vertex $v \in B$. By the choice of the indices, $(V_{s,t}^*,
V_{s',j}^*)$ is a $(d, \varepsilon)$-regular pair for all $j \in
[r]$. Therefore, there are at most $r\varepsilon |V_{s,t}^*|$
vertices in $V_{s,t}^*$ which have at most $(d -
\varepsilon)|V_{s',j}^*|$ neighbors in at least one of the sets
$V_{s',j}^*$ for $j \in [r]$. Since $v$ has at least $(1/3)m_{s,t}p$
neighbors in $V_{s,t}^*$, and
\begin{align} \label{eqn_eqnverticesB}
\frac{1}{3}m_{s,t}p  - r\varepsilon|V_{s,t}^*| \ge \frac{p}{3} |V_{s,t}^*| - r\varepsilon|V_{s,t}^*| > \frac{p}{4} |V_{s,t}^*|,
\end{align}
we can find one vertex $v_1 \in N_{G'}(v)$ which has at least $(d - \varepsilon)|V_{s',j}^*|$
neighbors in $V_{s',j}$ for all $j \in [r]$. We
will show by induction that there are $\Delta$ vertices $v_1, \ldots, v_\Delta$
which have many common neighbors in the sets $V_{s',j}^*$ for all $j
\in [r]$. Assume that for some $k \le \Delta -1$, we have found $v_1, \ldots, v_k \in N_{G'}(v)$ which have at
least $(d - \varepsilon)^k|V_{s',j}^*|$ common neighbors in
$V_{s',j}^*$ for all $j \in [r]$. For a fixed $j$, since
$(d - \varepsilon)^k|V_{s',j}^*| \ge \varepsilon|V_{s',j}^*|$ (recall that
we chose $\varepsilon \le (d/8)^\Delta$), by the $\varepsilon$-regularity
of pairs, there are at most $\varepsilon |V_{s,t}^*|$
vertices in $V_{s,t}^*$ which have less than $(d -
\varepsilon)^{k+1}|V_{s',j}^*|$ neighbors in $V_{s',j}^* \cap \bigcap_{i=1}^k
N_{G'}(v_i)$. Thus when we consider all the indices, there would be at most $r \varepsilon
|V_{s,t}^*|$ such `bad' vertices. By (\ref{eqn_eqnverticesB}),
since $(p/4)|V_{s,t}^*| - k > 0$, we can pick a vertex $v_{k+1} \in N_{G'}(v)$ not
equal to $v_1, \ldots, v_k$ so that the size of the common
neighborhood of $v_1, \ldots, v_{k+1}$ in $V_{s', j}^*$ is at least
$(d - \varepsilon)^{k+1}|V_{s',j}^*|$ for all $j \in [r]$. In the end,
we will find $v_1, \ldots, v_{\Delta}$ as promised.

Arbitrarily embed the neighbors of $h_v$ into $v_i$ one by one.
Since $H$ has maximum degree at most $\Delta$ and the neighborhood
of $h_v$ is an independent set, this embedding is a graph
homomorphism (note that we heavily rely on the fact that $N_H(h_v)$
is an independent set). Repeat it for other vertices of $B$. Since
$B$ is a set of constant size, and in (\ref{eqn_eqnverticesB}) we
have $(p/4)|V_{s,t}^*| - \Delta |B| > 0$, this can be done for every
vertex in $B$ even if they share the same set $V_{s,t}^*$. Moreover,
for two vertices $v,v' \in B$, their preimages $h_v = g^{-1}(v)$ and
$h_{v'} = g^{-1}(v')$ were chosen to be at distance at least 5 apart
from each other. Thus there will be no edges between the
neighborhood of $h_v$ and the neighborhood of $h_{v'}$.
Consequently, once we find a map as above for all the vertices in
the neighborhood of $W_B := g^{-1}(B)$, it will become a graph
homomorphism of $H[W_B \cup N(W_B)]$ to $G'$ (in fact, here we only
need $h^{-1}(v)$ and $h^{-1}(v')$ to be at distance 4 apart). Let $Z
= W_B \cup N(W_B)$.

For $v$ and $h_v$ as above, pick a vertex $w \in N_H^{(2)}(h_v)$. By
the fact $h_v \in \big( \cup_{j=1}^{r}W_{s',j} \big) \backslash
\big( \cup_{l=0}^{3} N^{(l)}(X) \big)$ and the property of the set
$X$ saying that edges not incident to $X$ only lie on $K_k^r$ in the
homomorphic image of $H$ into $R$, we know that $w \in W_{s',t'}$
for some $t' \in [r]$. Therefore by the condition on the size of the
common neighbors that we imposed on the images of $N_H(h_v)$, there
exists a set $C_w$ of size at least $(d - \varepsilon)^\Delta
|V_{s',t'}^*| \ge 4c |V_{s',t'}^*|$ inside $V_{s', t'}^*$ whose
every element is a possible image of $w$. Here we rely on the fact
the that vertices in $W_B$ are at distance at least 5 apart from
each other, since this implies that all the neighbors of $w$ in $Z$
are solely contained in $N_H(h_v)$, and thus all the vertices in
$C_w$ are indeed possible images of $w$. Even if we discard the
elements of $g(Z)$ from $C_w$, since $|Z| \le (\Delta+1)|B|$ is a
constant and $|V_{s',t'}^*|$ is linear in $n$, the size of the set
$C_w$ will be at least $2c |V_{s',t'}^*|$.
\end{pfof}

Equipped with this theorem, we can prove an embedding result for
general graphs $H$ which does not satisfy the condition of having
enough vertices with independent neighborhood. The following
corollary states that as long as the order of $H$ is slightly
smaller than that of $G$, we can still find a copy of $H$ in
subgraphs of $G(n,p)$. The necessity of $H$ being smaller than
$G(n,p)$ will be discussed in the next section.

\begin{cor} \label{cor_mainthmvar1}
For all integers $r, \Delta$, and reals $0 < p \leq 1$, $\gamma >
0$, there exists a constant $\beta > 0$ such that the following
holds. Let $H$ be an $r$-chromatic graph on at most
$n - 1/ \beta^2$ vertices with $\Delta(H) \leq \Delta$
and bandwidth at most $\beta n$. Then $G=G(n,p)$ a.a.s. satisfies
the following. Let $G' \subset G$ be a spanning subgraph with $\delta(G')
\geq (1 - 1/r + \gamma)np$, then $G'$ contains a copy of $H$.
\end{cor}
\begin{pf}
Let $\beta'=\beta_{\ref{thm_mainthm}}$ and $\beta = \beta'/2$.
Assume that $H$ is a graph with exactly $n -
1/\beta^2$ vertices which satisfies the condition above and label
the vertices as $1,\ldots, n - 1/\beta^2$ so that the bandwidth is
at most $\beta n$. We will construct a new graph $H'$ containing $H$
which satisfies the condition of Theorem \ref{thm_mainthm} with
parameter $\beta'$ as following. Insert an isolated vertex at the
end of every interval $[(\beta^2 n - 1)k + 1, (\beta^2 n -
1)(k+1)]$. Clearly, $H$ is still $r$-chromatic, since we added an
independent set. Moreover, since we added at most $1 / \beta^2$ new
vertices, $H'$ has at most $n$ vertices and bandwidth at most $\beta
n + 1/\beta^2 \leq \beta' n$. By the fact that all the new vertices
are isolated, for every $[a, a + \beta^2 n] \subset [1,n]$, there
exists a vertex with independent neighborhood. Since $\beta' \ge
\beta$, this also holds with $\beta$ replaced by $\beta'$.
Therefore we can apply Theorem \ref{thm_mainthmintro1} to find a
copy of $H'$ in $G$ which also gives us a copy of $H$ in $G$.
\end{pf}

\section{Packing Problem}
\label{section_applicationpackingproblem}

Throughout this section let $H_0$ be a fixed graph on $h$ vertices
with chromatic number $r$. We will investigate the following
problem: ``For a fixed $0 < p \leq 1$, when does every spanning $G'
\subset G(n,p)$ with $\delta(G') \geq (1 - 1/r + \gamma)np$ a.a.s.
have a perfect $H_0$-packing?''. Our goal is to extend the results
of Alon and Yuster \cite{MR1376050}, Koml{\'o}s, S{\'a}rk\"{o}zy,
and Szemer{\'e}di \cite{MR1829855} to random graphs. It is clear 
that $n$ must be a multiple of $h$ but is there any additional 
necessary condition?

For the simplicity of later exposition, before proceeding further,
we establish several properties of random graphs that will be used
later.

\begin{lemma} \label{lemma_randomgraphproperties}
Let $0 < p \le 1$ be fixed, and $C, \alpha$ be positive constants.
Then $G=G(n,p)$ satisfies the following properties with probability
$1 - e^{-\Omega_{\alpha, C, p}(n)}$.
\begin{enumerate}[(i)]
    \setlength{\itemsep}{1pt}
  \setlength{\parskip}{0pt}
  \setlength{\parsep}{0pt}
  \item Every vertex $v$ has degree $\deg(v) \in [(1-\alpha)np, (1+\alpha)np]$.
  \item Every pair of distinct vertices $v,w \in V$ have between $(1-\alpha)np^2$ and $(1+\alpha)np^2$ common neighbors.
  \item For all $X,Y \subset V$ of size $|X|,|Y| = \Omega(n), e(X,Y) \in [(1-\alpha)|X||Y|p, (1+\alpha)|X||Y|p]$. In particular,
  $e(X) = e(X,X)/2 = [(1 - \alpha)|X|^2p/2, (1 +
  \alpha)|X|^2p/2]$.
  \item For every set $X$ of size $|X| \le Cp^{-2}$,
  there are at most $e^{-\Omega_{\alpha}(|X|p)}n$ vertices $v
  \in V \setminus X$ which have $\deg(v, X) \notin [(1 -
  \alpha)|X|p, (1+\alpha)|X|p]$.
  \item For every set $X$ of size $|X| \le Cp^{-2}$,
  there are at most $e^{-\Omega_\alpha(|X|p^2)}n^2p$ edges $\{v,
  w\}$ in $G[V \setminus X]$ such that $v$ and $w$ have fewer
  than $(1 - \alpha)|X|p^2$ common neighbors in $X$.
\end{enumerate}
\end{lemma}
\begin{pf}
$(i),(ii),(iii)$ follows directly from Chernoff inequality and
taking union bounds. We omit the details. Let $X$ be a fixed set of
size $|X| \leq Cp^{-2}$. To prove $(iv)$, note that by Chernoff
inequality, the probability of a single vertex $v \in V \setminus X$
having $\deg(v, X) \notin [(1 - \alpha)|X|p, (1+\alpha)|X|p]$ is
$e^{-\Omega_\alpha(|X|p)}$. Thus the expected number of such
vertices in $V \setminus X$ is $e^{-\Omega_{\alpha}(|X|p)}n$. Since
these events for different vertices are mutually independent of each
other, we can apply Chernoff inequality once more to conclude that
with probability $1 - e^{-\Omega_{\alpha, C, p}(n)}$, there are at
most $2e^{-\Omega_\alpha(|X|p)}$ vertices $v \in V \setminus X$
which has $\deg(v, A) \notin [(1 - \alpha)|X|p, (1+\alpha)|X|p]$.
And since there are at most $\sum_{k \le Cp^{-2}} {n \choose k}$
choices for $X$, we can take the union bound to derive the
conclusion for all choices of $X$.

To prove $(v)$, first expose the edges between $X$ and $V \setminus
X$ and call a pair of vertices $\{v,w\} \in V \setminus X$
\textit{bad} if $v$ and $w$ have fewer than $(1 - \alpha)|X|p^2$
common neighbors in $X$. We will bound the number of bad pair of
vertices by bounding the number of pairs $\{v,w\}$ where (a) $v$ has
too few neighbors in $X$ or (b) $v$ has enough neighbors but $w$
does not have enough common neighbors with $v$ in $X$. To bound (a),
by $(iv)$ with $\alpha/2$ instead of $\alpha$, we know that there
are at most $e^{-\Omega_\alpha(|X|p)}n$ vertices $v \in V \setminus
X$ which have less than $(1 - \alpha/2)|X|p$ neighbors in $X$. Even
if we assume that all the pairs which contain these vertices are
bad, there will be at most $e^{-\Omega_\alpha(|X|p)}n^2$ such pairs.
Then to bound (b), assume that $v \in V \setminus X$ has more than
$(1 - \alpha/2)|X|p$ neighbors in $X$. Then by Chernoff inequality,
any $w \in V \setminus (X \cup \{v\})$ has at least $(1 -
\alpha/3)|N(v,X)|p \geq (1 - \alpha)|X|p^2$ neighbors in $X$ with
probability $1 - e^{-\Omega_\alpha(|X|p^2)}$. Since for distinct
vertices in $V \setminus (X \cup \{v\})$ these events are
independent, by using Chernoff inequality again, with probability $1
- e^{-\Omega_{\alpha, C, p}(n)}$, there will be at most
$e^{-\Omega_\alpha(|X|p^2)}n$ vertices $w \in V \setminus (X \cup
\{v\})$ such that $v$ and $w$ have fewer than $(1 - \alpha)|X|p^2$
common neighbors in $X$. By taking the union bound over all vertices
$v \in V \setminus X$, we can conclude that with probability $1 -
e^{-\Omega_{\alpha, C, p}(n)}$, the contribution from (b) is
$e^{-\Omega_\alpha(|X|p^2)}n^2$. Thus there are at most
$e^{-\Omega_\alpha(|X|p^2)}n^2$ bad pairs in $V \setminus X$.

Now expose the edges within $V \setminus X$. By Chernoff inequality,
with probability $1 - e^{-\Omega_{\alpha, C, p}(n^2p)}$, at most
$e^{-\Omega_\alpha(|X|p^2)}n^2p$ bad pairs will form an edge. Since
$n^2p \gg n$, all the required events happen with probability $1 -
e^{-\Omega_{\alpha, C, p}(n)}$. Since there are at most $\sum_{k\le
Cp^{-2}} {n \choose k}$ choices for $X$, we can take the union bound
to derive the conclusion for all choices of $X$.
\end{pf}

Coming back to our main question of this section regarding perfect
packing in subgraphs of $G(n,p)$, a simple observation combined with
Theorem \ref{thm_mainthmintro1} shows that $n$ being a multiple of
$h$ is sufficient for certain graphs. More precisely, this
condition is sufficient if $H_0$ contains a vertex having
independent neighborhood. To see this, let $n = h n'$ for some
integer $n'$ and let $H$ be a graph consisting of $n'$ vertex
disjoint copies of $H_0$. Then $H$ has bandwidth at most $h$ and
chromatic number $r$. Moreover, since each copy of $H_0$ has a
vertex with independent neighborhood, it is clear that the
conditions of Theorem \ref{thm_mainthmintro1} holds. Therefore
a.a.s.~$H$ can be embedded into every spanning subgraph $G' \subset
G(n,p)$ with $\delta(G') \geq (1 - 1/r + \gamma)np$. This result can
be formally stated as following.

\begin{prop} \label{prop_perfectpackingindepnbd}
Let $H_0$ be an $r$-chromatic graph which has a vertex not contained
in a triangle. Then for any $0 < p \leq 1$ and $\gamma > 0$, a.a.s.~
any spanning subgraph $G' \subset G(n,p)$ with minimum degree $(1 -
1/r + \gamma)np$ contains a perfect $H_0$-packing.
\end{prop}

In particular, if $H_0$ is a bipartite graph then a.a.s.~every $G'$
contains a perfect $H_0$-packing. One might suspect that the same
result holds for every graph $H_0$, but unfortunately this is not
true. The following proposition shows that perfect packing is
impossible for every graph which does not satisfy the condition of
having a vertex with independent neighborhood.

\begin{prop} \label{prop_noperfectpacking}
Let $H_0$ be a fixed graph whose every vertex is contained in a
triangle. Then for all $\varepsilon > 0$, there exists
$p_\varepsilon$ such that for all $0 < p \le p_\varepsilon$,
$G=G(n,p)$ a.a.s.~has a spanning subgraph $G'$ with $\delta(G') > (1
- \varepsilon)np$ such that at least $\varepsilon p^{-2}/3$ vertices
of $G'$ are not contained in a copy of $H_0$.
\end{prop}
\begin{pf}
Let $X$ be a set of size $|X| = \varepsilon p^{-2}/3$ and delete all
the edges of $G$ inside $X$. Since $X$ is a set of constant size,
the effect of these edges is asymptotically negligible. For a vertex
$v \in V \backslash X$, we expect that it has $|X|p = \varepsilon
p^{-1}/3$ neighbors in $X$, and by Lemma
\ref{lemma_randomgraphproperties} $(iv)$ with $\alpha=1$,
a.a.s.~there are at most $e^{-\Omega(\varepsilon p^{-1})}n$ vertices
in $V \backslash X$ which have degree greater than $2|X|p =
2\varepsilon p^{-1}/3$ into $X$. Let $W$ be the collection of all
such vertices and remove all the edges between $X$ and $W$. Note
that if $p \leq p_\varepsilon := c \varepsilon /
\log(\varepsilon^{-1})$ for sufficiently small constant $c$, then we
have $e^{-\Omega(\varepsilon p^{-1})}n \leq \varepsilon np/2$. Thus
we will not remove too many edges from any of the vertices in $X$
(and also from vertices in $W$ since $X$ is a set of constant size).

Then for all $y \in V \backslash (X \cup W)$ delete edges according
to the following rule. For every triangle $xyz$ in $G$ with $x \in
X$, $y,z \in V \setminus (X \cup W)$, remove the edge $yz$. By Lemma
\ref{lemma_randomgraphproperties} $(ii)$, a.a.s.~$x$ and $y$ have at
most $(9/8)np^2$ common neighbors. Moreover $\deg(y, X) \leq
2\varepsilon p^{-1}/3$ because $y \notin W$, and therefore we have
deleted at most $\big(2\varepsilon p^{-1}/3\big)\cdot(9/8)np^2 \leq
3\varepsilon np/4$ edges from $y$. Also note that there are no
further edges removed from $y$ since the process is symmetric and
$xyz$ forms a triangle if and only if $xzy$ forms a triangle. Let
$G'$ be the new graph. Then $\delta(G') \geq (1-\varepsilon)np$ and
the deleting process guarantees that every vertex $x \in X$ is not
contained in a triangle. However, since every point of $H_0$ is
contained in a triangle, there cannot exist a copy of $H_0$ in $G'$
which contains a vertex from $X$. Thus we have found a required
$G'$.
\end{pf}

\begin{rem}
It is easy to see that in this lemma, the constant $p$ must
be sufficiently small. Indeed, if $p$ is close
to 1, then every subgraph of $G(n,p)$ with minimum degree $(1 - 1/r
+ \gamma)np$ in fact has minimum degree greater than $(1 - 1/r +
\gamma/2)n$ and thus Koml{\'o}s, S{\'a}rk\"{o}zy, and
Szemer{\'e}di's theorem \cite{MR1829855} shows that a perfect
packing does exist.
\end{rem}

The consequence of this proposition is quite interesting. Given
$H_0$ as in the proposition, if we take $H$ to be the graph
consisting of $n'$ vertex disjoint copies of $H_0$ and let $n = h
n'$, then Proposition \ref{prop_noperfectpacking} is equivalent to
saying that, for any $\gamma < 1/r$ and sufficiently small $p$, $H$ a.a.s. cannot be embedded into some $G' \subset
G(n,p)$ with $\delta(G') \geq (1 - 1/r + \gamma)np$. Note that such
graph $H$ satisfies all the assumptions of Theorem
\ref{thm_mainthmintro1} except the one requiring $H$ to have enough
vertices with independent neighborhood. Therefore, this proposition
indicates the necessity of this condition for the theorem.

The proof of Proposition \ref{prop_noperfectpacking} also shows
that, even for arbitrary graph $H$ which is not necessarily a
disjoint union of copies of a fixed graph, if every vertex of $H$ is
contained in a triangle, then $H$ cannot be embedded into $G'$. On
the other hand, as we have seen from Corollary
\ref{cor_mainthmvar1}, if $H$ is allowed to be slightly smaller than
$G$ (constant difference is enough), then we can embed $H$ into the
subgraph $G'$. However, even though the required gap between the
sizes of $G$ and $H$ is only of constant size, this constant might
be rather huge because it comes from the regularity lemma. This
suggests a very natural question of determining the correct order of
magnitude of this gap. In the remaining part of this section, we
will investigate this question in the case when $H$ is the union of
vertex disjoint copies of $H_0$.

Let $K_{t_1, \ldots, t_r}$ be the complete $r$-partite graph with
parts having size $t_1, \ldots, t_r$ respectively. Next lemma shows
that for certain graphs, the assertion of Proposition
\ref{prop_noperfectpacking} is essentially best possible.

\begin{lemma} \label{lemma_VerticesNotInACopyIntermediate}
Let $H_0$ be the complete $r$-partite graph $K_{1, m, \ldots, m}$.
Then there exists a constant $C = C(r)$ such that for all $0 < p \leq 1$,
there exists $\varepsilon = \varepsilon(r,p)$ such that $G=G(n,p)$
a.a.s.~has the following property. For every spanning subgraph $G' \subset G$ with
minimum degree $\delta(G') \geq (1 - 1/r)np$, and every set $T
\subset V(G')$ of size $|T| \leq \varepsilon n$, all but at most
$Cp^{-2}$ vertices of $V \backslash T$ are contained in a copy of
$H_0$ in $G'$ which does not intersect $T$.
\end{lemma}
\begin{pf}
Let $V = V(G)$ and $\varepsilon = \varepsilon(r,p)$, $C =
C(\varepsilon)$ are constants which we choose later. Given $G'$
and $T$ as above, let $X \subset V \setminus T$ be an arbitrary set of size
$Cp^{-2}$, and let $Y = V \setminus (X \cup T)$. By assuming that the
events of Lemma \ref{lemma_randomgraphproperties} hold, we will show that
there exists a copy of $H_0$ in $G'$ which intersects $X$ but not
$T$.

For a vertex $x \in X$, let $N_x$ be the set of neighbors of $x$ in
$Y$ in the graph $G$, that is $N_x := N_G(x) \cap Y$, and note that
the size of $N_x$ is at least $(1 - 3\varepsilon p^{-1})np$ by Lemma
\ref{lemma_randomgraphproperties} $(i)$ and the fact $|X \cup T|
\leq 2 \varepsilon n$. Then in the graph $G'$, since the degree of
$x$ is at least $(1 - 1/r)np$, we can arbitrarily fix a set $N_x'
\subset N_{G'}(x) \cap Y$ of size $|N_x'| = (1 - 1/r - 2\varepsilon
p^{-1})np$. We claim that there exists a vertex $x \in X$ such that
$e_{G'}(N_x') \geq (1 - 1/(r-2) + \gamma)|N_x'|^2p/2$ for some
constant $\gamma > 0$. Then, by Corollary
\ref{cor_randomturanstrong}, $N_x'$ contains the complete
$(r-1)$-partite graph with parts of size $m$, which together with
$x$ will form a copy of $K_{1,m,\ldots,m}$ that intersects $X$ but
not $T$.

Thus it remains to verify the claim. To prove this claim we count
the number of triangles $xy_1y_2$ in $G'$ such that $x \in X,
y_1,y_2 \in N_x'$. Let this number be $M$. To lower bound $M$, first
bound the number of triangles $xy_1y_2$ in $G$ such that $x \in X,
y_1, y_2 \in N_x$, and $y_1y_2$ is an edge of the graph $G'$ (we will later subtract
the triangles whose $y_1$ or $y_2$ is not in $N_x'$).
Let this number be $M_0$. Since $|X
\cup T| \leq 2\varepsilon n$, by Lemma
\ref{lemma_randomgraphproperties} $(iii)$,
\[ e_{G'}(Y) \geq e_{G'}(V) - e_{G'}(V, X \cup T) \geq \left(1 - \frac{1}{r}\right)\frac{n^2p}{2} - e_{G}(V, X \cup T) \geq \left(1 - \frac{1}{r} - O(\varepsilon  p^{-1})\right)\frac{n^2p}{2}.\]
Let $\varepsilon' = \varepsilon'(r)$ be a small constant. If $C=
C(r)$ is large enough, by Lemma \ref{lemma_randomgraphproperties}
$(v)$, there are at most
$e^{-\Omega_{\varepsilon'}(|X|p^2)}n^2p=e^{-\Omega_{\varepsilon'}(C)}n^2p
= O(\varepsilon' n^2p)$ edges $\{v, w\}$ in $G[Y]$ which form a
triangle with fewer than $(1 - \varepsilon')C$ vertices $x \in X$.
These two facts provide the following bound on $M_0$:
\[ M_0 \ge \Big(e_{G'}(Y) - O(\varepsilon' n^2p) \Big) (1 - \varepsilon')C \ge \left(1 - \frac{1}{r} - O(\varepsilon p^{-1}) - O(\varepsilon') \right)\frac{Cn^2p}{2}. \]
To obtain a bound on $M$ from $M_0$, we can subtract the number of
triangles $xy_1y_2$ as above such that either $y_1$ or $y_2$ is
not in $N_x'$. Since $|N_x| = (1 - O(\varepsilon p^{-1}))np$,
\[ |N_x \setminus N_x'| = |N_x| - |N_x'| = \big(1 - O(\varepsilon p^{-1})\big)np - (1 - 1/r - 2\varepsilon p^{-1} )np = \big(1/r - O(\varepsilon p^{-1})\big)np. \]
Thus, if $\varepsilon = \varepsilon(p)$ is small enough, by Lemma \ref{lemma_randomgraphproperties} $(iii)$ we have,
\begin{align*}
 M &\ge M_0 - \sum_{x \in X} \left( e_{G'}(N_x \setminus N_x', N_x') + e_{G'}(N_x \setminus N_x')\right) \\
    &\ge M_0 - \sum_{x \in X} \left(1+O(\varepsilon p^{-1})\right)\left( \left(\frac{1}{r} - O(\varepsilon p^{-1})\right)\left( 1 - \frac{1}{r} + O(\varepsilon p^{-1}) \right)n^2p^3 + \left(\frac{1}{r} - O(\varepsilon p^{-1})\right)^2\frac{n^2p^3}{2} \right) \\
   &\ge \left(1 - \frac{1}{r} - O(\varepsilon p^{-1}) - O(\varepsilon') \right)\frac{Cn^2p}{2} - \sum_{x \in X} \left( \frac{1}{r}\left( 1 - \frac{1}{r} \right)n^2p^3 + \frac{1}{r^2}\frac{n^2p^3}{2} + O(\varepsilon p^{-1})n^2p^3 \right) \\
   &= \left(1 - \frac{3}{r} + \frac{1}{r^2} - O(\varepsilon p^{-1}) - O(\varepsilon') \right)\frac{Cn^2p}{2}.
\end{align*}

On the other hand we have, $M = \sum_{x \in X} e_{G'}(N_x')$. Thus
combining these two equations and using the fact $|X| = Cp^{-2}$, we
can find a vertex $x_0 \in X$ such that
\begin{align*}
 e_{G'}(N_{x_0}') \geq \frac{M}{|X|} &\geq \left(1 - \frac{3}{r} + \frac{1}{r^2} - O(\varepsilon p^{-1}) - O(\varepsilon') \right)\frac{n^2p^3}{2} \\
    &\geq \left( 1 - \frac{1}{r-2} + \gamma \right)\left( 1 - \frac{1}{r} \right)^2 \frac{n^2p^3}{2} \ge \left( 1 - \frac{1}{r-2} + \gamma \right)\frac{|N_{x_0}'|^2p}{2},
\end{align*}
for some constant $\gamma > 0$ depending on $r$, small enough
$\varepsilon'$ depending on $r$, and
$\varepsilon$ depending on $r$ and $p$. This concludes the proof.
\end{pf}

Next, we extend Lemma \ref{lemma_VerticesNotInACopyIntermediate} to all graphs $H_0$.

\begin{lemma} \label{lemma_VerticesNotInACopy}
Let $H_0$ be a fixed $r$-chromatic graph. Then there exists
a constant $C=C(r)$ such that for every $0 < p \leq
1$, there exists $\varepsilon = \varepsilon(r, p)$ such that $G=G(n,p)$ a.a.s.
has the following property. For every spanning subgraph $G' \subset
G$ with minimum degree $\delta(G') \geq (1 - 1/r)np$, and every set
$T \subset V(G')$ of size $|T| \leq \varepsilon n$, all but at most
$Cp^{-2}$ vertices of $V \backslash T$ are contained in a copy of
$H_0$ in $G'$ which does not intersect $T$.
\end{lemma}
\begin{pf}
Let $V = V(G)$, and
$C=C_{\ref{lemma_VerticesNotInACopyIntermediate}}(r)$. Let
$\varepsilon \leq
\varepsilon_{\ref{lemma_VerticesNotInACopyIntermediate}}(r,p)$ and
$D = D(r,p,\varepsilon)$ be constants to be chosen later. We may
assume that $H_0$ is a complete $r$-partite graph with equal parts
of size $s$. Throughout the proof we condition on the event that the
statements of Lemma \ref{lemma_randomgraphproperties} holds.

Given $G'$ and $T$ as above, let $X \subset V \setminus T$ be an arbitrary set
of size $Cp^{-2}$. We will show that there exists a copy of $H_0$ in
$G'$ which intersects $X$ but not $T$. By Lemma
\ref{lemma_VerticesNotInACopyIntermediate} we can find a complete
$r$-partite graph with parts $\{x\} \cup Z_1 \ldots \cup Z_{r-1}$
such that $x \in X$ and $|Z_i| = Dsp^{-1}$ (note that in Lemma
\ref{lemma_VerticesNotInACopyIntermediate}, the part size $m$ can be
an arbitrary constant). Let $Z = Z_1 \cup \ldots \cup Z_{r-1}$ and
$Y = V \backslash (X \cup Z \cup T)$. Note that $|Y| \geq (1 -
2\varepsilon)n$ for large enough $n$. We construct a set $A \subset Y$ of
size $s-1$ and sets $Z_i' \subset Z_i$ of size $s$ for $1 \leq i
\leq r-1$ such that $A \cup Z_1' \cup \ldots \cup Z_{r-1}'$ forms a
complete $r$-partite graph.

By Lemma \ref{lemma_randomgraphproperties} $(iv)$, there are at most
$e^{-\Omega_\varepsilon(Ds)}n$ vertices in $V \setminus X$ such that
$\deg_{G}(y, Z_i) > (1 + \varepsilon)|Z_i|p = (1+\varepsilon)Ds$ for
any fixed $1 \leq i \leq r-1$. Hence if $D = D(\varepsilon, p)$ is
large enough, there are at most $re^{-\Omega_\varepsilon(Ds)}n =
O(\varepsilon n p)$ vertices $y \in Y$ which has $\deg_{G}(y, Z_i) >
(1+\varepsilon)Ds$ for at least one $1 \leq i \leq r-1$. Let $Y_0$
be these vertices. Then we have the crude bound $e_{G'}(Y_0, Z) \le
O(\varepsilon np)|Z|$. Let $Y_1$ be the collection of vertices in $Y
\setminus Y_0$ which have at least $\varepsilon |Z_i|p =
\frac{\varepsilon|Z|p}{r-1}$ neighbors in $Z_i$ in the graph $G'$
for all $1 \leq i \leq r-1$, and $Y_2 := Y \setminus (Y_0 \cup
Y_1)$. Then since $Y_1 \subset Y \setminus Y_0$,
\[ e_{G'}(Y_1, Z) \le \sum_{i=1}^{r-1} e_{G'}(Y_1, Z_i) \le \sum_{i=1}^{r-1} (1+\varepsilon)|Y_1||Z_i|p = |Y_1| \cdot (1+\varepsilon)|Z|p, \]
and since $Y_2 = Y \setminus (Y_0 \cup Y_1)$,
\[ e_{G'}(Y_2, Z) \le |Y_2| \cdot \left((1 +
\varepsilon)\frac{|Z|p}{r-1}(r-2) + \varepsilon
\frac{|Z|p}{r-1}\right). \]
Thus we have,
\begin{align*}
e_{G'}(Y, Z) &\leq e_{G'}(Y_0, Z) + e_{G'}(Y_1, Z) + e_{G'}(Y_2, Z) \\
&\leq O(\varepsilon np)|Z| + |Y_1| \cdot (1+\varepsilon)|Z|p + n \cdot \left((1 +
\varepsilon)\frac{|Z|p}{r-1}(r-2) + \varepsilon
\frac{|Z|p}{r-1}\right) \\
             &= \left( \frac{|Y_1|}{n} + \frac{r-2}{r-1} + O(\varepsilon) \right)|Z|np.
\end{align*}
On the other hand, by the minimum degree condition of $G'$,
\begin{align*}
e_{G'}(Y, Z) &=  \sum_{z \in Z} \left( \deg_{G'}(z,V) - \deg_{G'}(z,V \setminus Y) \right) \\
&\geq \left( \left(1 - \frac{1}{r} \right)np - (n - |Y|)\right)|Z|
    \geq \left(1 - \frac{1}{r} - 2 \frac{\varepsilon}{p} \right)|Z|np. \nonumber
\end{align*}
By combining the previous inequalities and dividing each side by
$|Z|np$ we have,
\begin{align*}
\frac{|Y_1|}{n}  \geq 1 - \frac{1}{r} - 2 \frac{\varepsilon}{p} - \frac{r-2}{r-1} - O(\varepsilon) \geq \frac{1}{2r(r-1)}.
\end{align*}
The last inequality holds if we pick $\varepsilon=\varepsilon(r,p)$
small enough. Thus there are at least $\frac{1}{2r(r-1)}n$ vertices which have at least $\varepsilon
|Z_i|p = \varepsilon Ds$ neighbors in $Z_i$ for all $1 \leq i \leq
r-1$. Let $D \geq \varepsilon^{-1}$, and for each such vertex fix
$s$ points in each $Z_i$ which are adjacent to that vertex. Since
there are only $\binom{|Z_i|}{s}$ possible subsets of size $s$ in
each $Z_i$, and these numbers are constants, if $n$ is large enough
then by pigeonhole principle we can find $s-1$ vertices $y_1, y_2,
\ldots, y_{s-1}$ which are adjacent to the same $s$-tuple of
vertices in every $Z_i$. Let $A = \{ y_1, y_2, \ldots, y_{s-1} \}$,
and for each $i$, let $Z_i'$ be the $s$-tuple which is adjacent to
these vertices. Recall that $x \in X$ was a vertex chosen at the beginning,
which forms a complete $r$-partite graph together with the sets
$Z_1, \ldots, Z_{r-1}$. Since $Z_i'$ are subsets of $Z_i$,
 $Z_1' \cup \ldots Z_{r-1}' \cup \big( A \cup
\{x\} \big)$ forms a complete $r$-partite graph with $s$ vertices in
each parts which intersects $X$ but not $T$.
\end{pf}

We are now ready to prove the following theorem which when combined with
Proposition \ref{prop_noperfectpacking} establishes Theorem
\ref{thm_packingtheorem}.

\begin{thm} \label{thm_correctpacking}
Let $H_0$ be a fixed $r$-chromatic graph. There exists a constant $C
= C(r)$ such that for every fixed $0 < p \leq 1$ and $\gamma>0$, if a
spanning subgraph $G'\subset G$ satisfies $\delta(G')\geq
(1-1/r+\gamma)np$, then $G'$ contains vertex disjoint copies of
$H_0$ covering all but at most $Cp^{-2}$ vertices.
\end{thm}
\begin{pf}
Let $\Delta = \Delta(H_0)$ and $h = |V(H_0)|$. Let
$C=\max\{2C_{\ref{lemma_VerticesNotInACopy}}, 2rh\}$,
$d=d_{\ref{lemma_lemmaforG}}(r,p,\gamma)$, and $b_0 =
b_{\ref{lemma_lemmaforG}}$. Then let
\[ \varepsilon = \frac{1}{2}\min \left\{\varepsilon_{\ref{thm_blowuplemma}}(\frac{d}{2},\Delta,c,r),  \varepsilon_{\ref{lemma_lemmaforG}}(r,p,\gamma), \varepsilon_{\ref{lemma_VerticesNotInACopy}}(r,p), \frac{d}{2} \right\}, \]
and $\xi = \xi_{\ref{lemma_lemmaforG}}(r,p,\gamma,\varepsilon)$.

Assume that $G' \subset G(n,p)$ is given as above. Lemma
\ref{lemma_lemmaforG} applied to $G'$ provides us a subgraph $G''
\subset G'$, a graph $R$ over the vertex set $[k] \times [r]$, a set
$B$ with $|B| = b \leq b_0$, sets $(V_{i,j}^*)$, and a $r$-equitable
integer partition $(m_{i,j})_{1 \leq i \leq k, 1 \leq j \leq r}$
satisfying $(i), (ii), (iii), (iv)$ of Lemma \ref{lemma_lemmaforG}. Let $n' := n - |B| = \sum_{i,j}
m_{i,j}$ be the number of vertices not in $B$. Since copies of $H_0$
in $G''$ are also copies in $G'$, by abusing notation, we denote
$G'$ for the graph $G''$. Note that by doing this, we can only
guarantee $\delta(G') \ge (1 - 1/r + 4\gamma/5)np$.

We first find copies of $H_0$ containing vertices of $B$ and only
using vertices from $B \cup \Big(\bigcup_{i,j} V_{i,j}^* \Big)$. Let
$T = V \backslash \left(B \cup \big( \cup_{1 \leq i \leq k, 1 \leq j
\leq r} V_{i,j}^* \big)\right)$ and note that
\[ |T| \leq n - (1 - \varepsilon)\sum_{i,j}m_{i,j} - |B| \leq n - (1 - \varepsilon)(n - |B|) - |B|
\le \varepsilon n \leq \varepsilon_{\ref{lemma_VerticesNotInACopy}}
n. \] By Lemma \ref{lemma_VerticesNotInACopy} if $|B| \geq
(C/2)p^{-2}$ then we can find a copy of $H_0$ in $G'$ which
intersects $B$ but does not intersect $T$. Move the vertices of this
copy to $T$. Repeat this process, as long as $|B| \geq (C/2)p^{-2}$,
one can find a copy of $H_0$ intersecting $B$ but not $T$ (note that
$|T| \leq \varepsilon n + |B|h \leq
\varepsilon_{\ref{lemma_VerticesNotInACopy}} n$ at any point of this
process). In the end we will have vertex disjoint copies of $H_0$
and at most $(C/2)p^{-2}$ vertices left in $B$. The leftover
vertices of $B$ will remain uncovered. Our next task is to find a
$H_0$-packing in the remaining part. Let $S$ be the vertices
belonging to the copies of $H_0$ found so far.

Let $\delta_{i,j} = |V_{i,j}^* \cap S|$ and construct $(n_{i,j})_{1
\leq i \leq k, 1 \leq j \leq r}$ as following. For $i \in [k-1]$,
let $t_i$ be the largest integer smaller than $\min_{s} (m_{i,s} -
\delta_{i,s})$ which is divisible by $h$, and let $n_{i,j} = t_i$
for all $j \in [r]$. Then pick $t_k$ so that
$\sum_{i=1}^{k}\sum_{j=1}^{r} (t_i + \delta_{i,j}) \in (n' - rh,
n']$ is divisible by $rh$. Recall that $\sum_{i=1}^{k}\sum_{j=1}^{r}
m_{i,j} = n'$. Since $|m_{i,j} - m_{i,j'}| \le 1$ for all $i,j ,j'$,
we are modifying each $m_{i,j}$ by at most $\max_{s,t}(\delta_{s,t}
+ 1) + rh$ to construct $n_{i,j}$ for all $i,j$. Since $\delta_{i,j}
\le |S|$ for all $i,j$ and $S$ has constant size, it shows that
$|m_{i,j} - n_{i,j}|$ is at most some constant. Thus the integer
partition $(n_{i,j})$ satisfies the following properties.
\begin{enumerate}[(i)]
  \setlength{\itemsep}{1pt}
  \setlength{\parskip}{0pt}
  \setlength{\parsep}{0pt}
\item $n' - \sum_{i,j} \delta_{i,j} \ge \sum_{i,j} n_{i,j} \ge n' - \sum_{i,j} \delta_{i,j} - rh$,
\item $n_{i,j} \in [m_{i,j} -
(\xi/2)n, m_{i,j} + (\xi/2)n]$,
\item $n_{i,j} = n_{i,j'}$ for all
$1 \leq i \leq k, 1 \leq j, j' \leq r$, and
\item $h$ divides
$n_{i,j}$ for all $1 \leq i \leq k,  1\leq j \leq r$.
\end{enumerate}
It then follows that $n_{i,j} + \delta_{i,j} \in [m_{i,j} - \xi n, m_{i,j} + \xi
n]$. So by Lemma \ref{lemma_lemmaforG} we can find sets $V_{i,j}$
such that $|V_{i,j}| \geq n_{i,j} + \delta_{i,j}$ which are
$(d,\varepsilon)$-super-regular on $K_k^r$. Let $V_{i,j}' = V_{i,j}
\backslash S$, and we have $|V_{i,j}'| \geq |V_{i,j}| - \delta_{i,j}
\ge n_{i,j}$. Remove some vertices so that $|V_{i,j}'| = n_{i,j}$.
The number of removed vertices is at most $n' - \sum_{i,j}(n_{i,j}+\delta_{i,j}) \le rh$.
These vertices together with the remaining vertices of $B$
will form the $(C/2)p^{-2} + rh \le Cp^{-2}$
uncovered vertices. Further note that we removed only at most
some constant number of vertices from each $V_{i,j}$ to obtain $V_{i,j}'$.

Since $(V_{i,j})_{1 \leq i \leq k, 1 \leq j \leq r}$ is
$(d,\varepsilon)$-super-regular on $K_k^r$ and we removed only at most
constant number of vertices from each part to obtain $V_{i,j}'$, we can conclude that
$(V_{i,j}')_{i,j}$ is $(d-\varepsilon,2\varepsilon)$-super-regular
on $K_k^r$ (Lemma \ref{lemma_preservesuperregularity}).
Thus we may apply the blow-up lemma to the super-regular
partitions $(V_{i,j}')_{1 \leq j \leq r}$ for each fixed $i \in
[k]$ to find a perfect $H_0$-packing in each of them. By (iii) and
(iv) of the previous paragraph, it suffices to show that the
complete $r$-partite graph with $h$ vertices in each
class contains a perfect $H_0$-packing, or equivalently, $r$ vertex disjoint
copies of $H$ has an $r$-coloring in which every color class has size $h$.
Assume that $H_0$ has an $r$-coloring with color classes of size $h_1, \ldots, h_r$.
Then by renaming the colors, we can color the $i$-th copy of $H_0$ so that the $j$-th color
class of it has $h_{i+j-1}$ vertices (addition of indices are modulo $r$).
In this way, we will end up with a coloring of $r$ vertex disjoint copies
of $H_0$ in which every color class has size $\sum_{i=1}^r h_i = h$.
\end{pf}

\section{Concluding Remarks}
\label{section_concludingremark}

\noindent $\bullet$ \hspace{0.1cm} In this paper, we proved that for
all integers $r$ and $p \in (0, 1]$, there exists $\beta$ such that
if $H$ is an $r$-chromatic graph on $n$ vertices with bounded
degree, bandwidth at most $\beta n$, and has enough vertices whose
neighbors form an independent set, then $G(n,p)$ a.a.s. has the
following property. Every spanning subgraph $G' \subset G(n,p)$ with
minimum degree at least $(1 - 1/r + \gamma)np$ contains a copy of
$H$. It would be interesting to know whether this theorem holds for
$p \ll 1$ or not. As mentioned in the introduction, B\"{o}ttcher,
Kohayakawa, and Taraz \cite{MR000006} proved that for fixed $\eta,
\gamma > 0, \Delta > 1$ there exist positive constants $\beta$ and
$c$ such that if $p \geq c(\log n / n)^{1 / \Delta}$ then a.a.s
every subgraph of $G(n,p)$ with minimum degree at least $(1/2 +
\gamma)np$ contains a copy of any bipartite graph $H$ with
$(1-\eta)n$ vertices, maximum degree $\Delta$ and bandwidth at most
$\beta n$. However, it is plausible that one can even embed a
spanning bipartite graph $H$ under the same conditions. The
technique we used in this paper cannot be applied mainly because of
the lack of the corresponding blow-up lemma in the range $p \ll 1$.
It is hopeful that a sparse version of the blow-up lemma (if one
exists) will allow us to extend the same proof.

\medskip

\noindent $\bullet$ \hspace{0.1cm} In view of the results of Koml\'os
\cite{Komlos2}, Shokoufandeh and Zhao \cite{ShZh1}, \cite{ShZh2},
and K\"{u}hn and Osthus \cite{KuOs}, which establishes the best
possible minimum degree condition for packing problems, it is likely
that in Theorem \ref{thm_packingtheorem},
the minimum degree condition $(1-1/r+\gamma)np$ can be further
relaxed. However, we did not further pursue towards this direction as
our primary goal was to study the packing problem in
connection to Theorem \ref{thm_mainthmintro1}.

\medskip

\noindent $\bullet$ \hspace{0.1cm} For a graph $G$, let $\lambda_1
\ge \lambda_2 \ge \ldots \ge \lambda_n$ be the eigenvalues of its
adjacency matrix. The quantity $\lambda(G) = \max\{\lambda_2, -
\lambda_n\}$ is called the second eigenvalue of $G$. A graph $G =
(V,E)$ is called an $(n,d,\lambda)$-graph if it is $d$-regular, has
$n$ vertices, and the second eigenvalue of $G$ is at most $\lambda$.
It is well known (see e.g., survey \cite{MR2223394}) that if
$\lambda$ is much smaller than the degree $d$, then $G$ has certain
random-like properties. Thus $\lambda$ could serve as some kind of
``measure of randomness'' in $G$. By using an almost identical
argument as in the proof of Theorem \ref{thm_mainthmintro1} we can
prove the bandwidth theorem for pseudorandom graphs as well.

\begin{thm} \label{thm_pseudomainthm}
For all integers $r, \Delta$, and reals $\gamma > 0$ and $0 < p \leq
1$, there exists a constant $\beta > 0$ such that, for an
$(n,d,\lambda)$ graph $G$ with $d = np$ and $\lambda = o(n)$, if $n$
is large enough, then any spanning subgraph $G' \subset G$ with
minimum degree $\delta(G') \geq (1 - 1/r + \gamma)np$ contains a
copy of every graph $H$ on $n$ vertices which satisfies the
following properties. $(i)$ $H$ is $r$-chromatic, $(ii)$ has maximum
degree at most $\Delta$, $(iii)$ has bandwidth at most $\beta n$
with respect to a labeling of vertices by $1,2,\ldots, n$, and
$(iv)$ for every interval $[a, a + \beta^2 (n/\lambda)] \subset [1,
n]$, there exists a vertex $v \in H$ such that $N_H(v)$ is an
independent set.
\end{thm}

The sketch of the proof will be given in the appendix. \\

\noindent {\bf Acknowledgment.} We thank the anonymous referees for
their careful reading and valuable comments. We would also like to
thank Wojciech Samotij and Deryk Osthus for their helpful comments.

\appendix

\section{Pseudo-random Graphs}
\label{section_pseudorandomgraphs}

In this appendix, we give an outline of the proof of Theorem
\ref{thm_pseudomainthm}. We can use an argument identical to the
proof of Theorem \ref{thm_mainthmintro1} given in Section
\ref{section_maintheorem}. Note that in the proof of Theorem
\ref{thm_mainthmintro1}, all the properties of random graphs were
encoded in Lemma \ref{lemma_lemmaforG}, Lemma for $G$. Thus in order
to prove Theorem \ref{thm_pseudomainthm}, we need to prove the
following lemma which is a version of Lemma \ref{lemma_lemmaforG}
for pseudorandom graphs.

\begin{lemma} \label{lemma_lemmaforGpseudo}
For all integer $r$, and reals $0 < p \leq 1$, $\gamma > 0$, there
exists $d> 0$ and $\varepsilon_0 > 0$ such that for every positive
$\varepsilon \leq \varepsilon_0$ there exists $b_0$, $K_0$ and
$\xi_0 > 0$ such that any $(n,d,\lambda)$-graph $G$ with $d=np$ and
$\lambda = o(n)$ satisfies the following if $n$ is large enough. For
every subgraph $G' \subset G$ with $\delta(G') \geq ((r-1)/r +
\gamma)np$ there exist a subgraph $G'' \subset G'$ with $\delta(G'')
\geq (1-1/r + 4\gamma/5)np$, a set $B$ of size at most $b_0
\lambda$, an $r$-equitable integer partition $(m_{i,j})_{1 \leq i \leq k, 1 \leq j \leq r}$ of $n - |B|$, sets $( V_{i,j}^* )_{1 \leq i \leq k, 1 \leq j \leq r}$,
and a graph $R$ on vertex set $[k] \times [r]$ with $k \leq K_0$ such that,\\
\begin{tabular}{cl}
$(i)$ & $K_k^r \subset C_k^r \subset R$ and $\delta(R) \geq (1-1/r + \gamma/2)kr$, \\
$(ii)$ &  $\forall 1 \leq i \leq k, 1 \leq j \leq r$, $m_{i,j} \geq (1-\varepsilon)n/(kr)$,\\
$(iii)$ & $\forall 1 \leq i \leq k, 1 \leq j \leq r$, $m_{i,j} \ge |V_{i,j}^*| \geq (1-\varepsilon)m_{i,j}$,\\
$(iv)$ & $(V_{i,j}^*)_{1 \leq i \leq k, 1 \leq j \leq r}$ is $(d,\varepsilon)$-regular on $R$ in $G''$, such that
\end{tabular}

\noindent for every choice of $(n_{i,j})_{1 \leq i \leq k, 1 \leq j \leq r}$
 with $m_{i,j} - \xi_0 n \leq n_{i,j} \leq m_{i,j} + \xi_0 n$ and $\sum_{i,j} n_{i,j} \leq n - |B|$,
there exists a partition $(V_{i,j})_{1 \leq i \leq k, 1 \leq j \leq r}$ of $V
\backslash B$
with \\
\begin{tabular}{cl}
$(a)$ & $|V_{i,j}| \geq n_{i,j}$, $V_{i,j}^* \subset V_{i,j}$, $\forall 1 \leq i \leq k, 1 \leq j \leq r$, \\
$(b)$ & $(V_{i,j})_{1 \leq i \leq k, 1 \leq j \leq r}$ is $(d,\varepsilon)$-regular on $R$ in $G''$ and \\
$(c)$ & $(V_{i,j})_{1 \leq i \leq k, 1 \leq j \leq r}$ is $(d,\varepsilon)$-super-regular on $K_k^r$ in $G''$.
\end{tabular}
\end{lemma}

Note that the size of the set $B$ is now $b_0 \lambda$ compared to
a constant in the random graph case. This subtlety slightly affects the proof
of Theorem \ref{thm_pseudomainthm}, and that is why we now need a
vertex with independent neighborhood in every interval
of length $\beta^2(n/\lambda)$ instead of $\beta^2 n$ as in the random
graph case (see the statement of Theorem \ref{thm_pseudomainthm}).
The proof of Lemma \ref{lemma_lemmaforGpseudo} is almost identical
to the proof of Lemma \ref{lemma_lemmaforG}. However, in order to apply
the same argument, we need lemmas corresponding to Lemma \ref{lemma_reducedgraphmindegree} and
Lemma \ref{lemma_badset} for pseudorandom graphs.

The main fact
that we use about $(n,d,\lambda)$-graphs is the following formula
established by Alon (see, e.g., \cite{MR2223394}) which connects
between eigenvalues and edge distribution.
\begin{lemma} \label{lemma_pseudoedgedistribution}
If $G=(V,E)$ is an $(n,d,\lambda)$-graph, then for any $X,Y \subset V$ we have,
\[ \Big|e(X,Y) - \frac{d}{n}|X||Y|\Big| \leq \lambda \sqrt{|X||Y|}. \]
\end{lemma}

In particular, if $G$ is a $(n,d,\lambda)$-graph for $d=np$ and $\lambda = o(n)$ as in Theorem
\ref{thm_pseudomainthm}, then for linear size sets $|X|, |Y| \geq
\varepsilon n$, we have $e(X,Y) = (1+o(1))|X||Y|p$. The following
lemmas are a pseudorandom graph version of Lemma \ref{lemma_reducedgraphmindegree} and
Lemma \ref{lemma_badset}.

\begin{lemma} \label{lemma_reducedgraphmindegreepseudo}
Let $0< p \le 1$ and $\alpha, \gamma>0$ be fixed. There exists
$\epsilon_0$ such that for all
$\varepsilon \le \varepsilon_0$ and $d > 0$, the following holds.
Given a $(n,d,\lambda)$-graph $G$ with $d = np$ and $\lambda = o(n)$,
let $V_0, V_1, \ldots, V_k$ be a pure
$(d,\varepsilon)$-regular partition of a subgraph $G' \subset G$, and
$R$ be its reduced graph. For large enough $n$, if $G'$ has minimum degree at least
$(\alpha + \gamma)np$, then $R$ has minimum degree at least $(\alpha
+ 3\gamma/4)k$.
\end{lemma}
The proof of this lemma is a line by line translation of the proof
of Lemma \ref{lemma_reducedgraphmindegree} once we notice that
$e_G(V_i, V_j) \le (1 + \varepsilon)|V_i||V_j|$ for all $i,j \ge 1$ by
Lemma \ref{lemma_pseudoedgedistribution}.

\begin{lemma}\label{lemma_badsetpseudo}
Let $0 < p \leq 1$ be fixed and $T$ be an integer. Then for every
$\varepsilon>0$, there exists a constant $b_0=b_0(p,T,\varepsilon)$
such that any $(n,d,\lambda)$-graph $G$ with $d = np$ and $\lambda=o(n)$
satisfies the following. For arbitrary subsets $V_1, \ldots, V_T$ of
the vertex set $V$ with $|V_i| \geq \varepsilon n$ for all $1 \le i
\le T$, there exists a set $B$ of size at most $b_0 \lambda$ such
that for all $v \in V \setminus B$, we have $\deg(v, V_i) \in [(1 -
\varepsilon) |V_i|p, (1 + \varepsilon)|V_i|p]$ for all $1 \le i \le
T$.
\end{lemma}
\begin{pf}
Let $m := |V_i| \ge \varepsilon n$ and $b = 2/(\varepsilon^2 p)$. Fix an index $i$ and let $X_i$ be the
collection of vertices which
have degree greater than $(1 + \varepsilon)|V_i|p$ in $V_i$.
Assume that $|X_i| \ge b \lambda$. By discarding some of the vertices, we may assume that $|X_i| = b \lambda$.
Then by how we chose $X_i$, and Lemma \ref{lemma_pseudoedgedistribution}, and $\lambda = o(m)$,
\[ e(V_i, X_i) > (1 + \varepsilon)|X_i|mp - e(X_i, X_i) \ge (1 + \varepsilon)|X_i|mp - |X_i|^2p - \lambda |X_i| \ge (1 + \varepsilon/2)|X_i|mp. \]
On the other
hand, by Lemma \ref{lemma_pseudoedgedistribution}, we have,
\[ e(V_i, X_i) \leq p|V_i||X_i| + \lambda \sqrt{|V_i||X_i|} \leq |X_i|mp + \lambda n.  \]
Thus by combining these two bounds we get $(\varepsilon/2)|X_i|mp < \lambda n \le \lambda m / \varepsilon$.
Recall that $|X_i| = b \lambda$, and hence this give $b < 2/(\varepsilon^2 p)$ which
is a contradiction.
By the definition of $X_i$, this shows that there are at most $b \lambda$ vertices which
have degree greater than $(1 + \varepsilon)|V_i|p$ in $V_i$.
Same estimate holds for the vertices which
have degree less than $(1 - \varepsilon)|V_i|p$ in $V_i$ and for other indices as well.
Therefore, the total number of vertices which does not satisfy $\deg(v, V_i) \in [(1 -
\varepsilon) |V_i|p, (1 + \varepsilon)|V_i|p]$ for at least one
index $i$ is at most $2T b \lambda$. Let $b_0 = 2Tb$ and this concludes the proof.
\end{pf}

We omit the proof of Lemma \ref{lemma_lemmaforGpseudo} and the
deduction of Theorem \ref{thm_pseudomainthm} from it. To prove Lemma
\ref{lemma_lemmaforGpseudo}, one should notice that in the proof of
Lemma for G (Lemma \ref{lemma_lemmaforG}), once we apply Lemma
\ref{lemma_reducedgraphmindegree} and Lemma
\ref{lemma_badset}), we do not use any other property of random
graphs, and since we have the corresponding lemmas for pseudorandom
graphs, the exact same proof will work for pseudorandom graphs as
well. Then in order to prove Theorem \ref{thm_pseudomainthm}, we can
first apply lemma for $G$ as in the proof of Theorem
\ref{thm_mainthmintro1}. Since we need not develop a new lemma for
$H$, the remaining proof of Theorem \ref{thm_pseudomainthm} will be
the same as in the random case. One difference is that $B$ is now a
set of sublinear size compared to a set of constant size in the
random case. However, the only part where we actually needed $B$ to
be a set of constant size instead of sublinear size was where we
wanted vertices with independent neighborhood (Claim
\ref{clm_verticesforB}). This issue has been resolved by increasing
the number of such vertices accordingly in the graph $H$.

\end{document}